\documentclass[leqno,12pt]{amsart}
\usepackage{eepic,epic}
\usepackage{amsmath}

\theoremstyle{plain}
\newtheorem{lemma}{Lemma}[section]
\newtheorem{theorem}[lemma]{Theorem}
\newtheorem{assertion}[lemma]{Assertion}
\newtheorem*{corollary}{Corollary}

\theoremstyle{definition}
\newtheorem*{definition}{Definition}

\theoremstyle{remark}
\newtheorem{remark}{Remark}

\newtheorem*{note}{Note}

\def\exp{{\operatorname{exp}}}

\def\Z{{\mathbb Z}}

\def\R{{\mathbb R}}

\def\A{{\mathbb A}}
\def\G{{\mathbb G}}
\def\al{\alpha}

\begin{document}

\title[principal $\Gamma$-cone]
{Principal $\Gamma$-cone for a tree
}
\author{Kyoji Saito}
\address{Research Institute for Mathematical Sciences, Kyoto University, 
Kyoto 606-8502, Japan.}

\maketitle
 \vspace{-1.5cm}

\begin{abstract} 

For any tree $\Gamma$, we introduce $\Gamma$-cones and 
solve a problem of enumeration of chambers contained 
in two particular,\! 
called the {\it principal}, $\Gamma$\!-cones. 
The problem can be reformulated
combinatorially: the enumeration of linear extensions of 
the two bipartite (called principal) orderings on $\Gamma$.
We characterize the principal $\Gamma$\!-cones by the\! {\it strict maximality} of 
the number of their chambers, and 
give a formula for this maximal number by 
a {\it finite sum of hook length formulae}. 
We explain the formula through  
the {\it block decomposition} of the principal $\Gamma$-cones.
The results have their origin and applications in the study of topology 
related to Coxeter groups.


\end{abstract}


Let $V_\Pi\!:=\!\oplus_{\al\in\Pi} \R v_\al/\R\!\sum_{\al\in\Pi}\!v_\al$ 
be the space of all configurations in $\R$ 
labeled by a finite set $\Pi$ 
(see \S1),
on which the permutation group $\mathfrak{S}(\Pi)$ acts irreducibly 
as a reflection group. 
The space $V_\Pi$ decomposes into chambers cut by the system of 
all hyperplanes fixed by reflections on $V_\Pi$ associated to  
the transpositions in  $\mathfrak{S}(\Pi)$. 
Given a  graph $\Gamma$, having $\Pi$ 
as its vertex set, the system of reflection hyperplanes 
associated to the transpositions of vertices on 
the edges of $\Gamma$ cuts $V_\Pi$ into components.
Each component, which we call a {\it $\Gamma$-cone},
is subdivided into chambers. 
If
$\Gamma$ is a tree, through a use of the 
principal decomposition of $\Gamma$ (\S3), 
we choose two particular 
$\Gamma$-cones, called the {\it principal $\Gamma$-cones}.
The main results of this article are 
i) {\it a characterization of the principal cones 
by the strict maximality of the number of their chambers} 
(Theorem 3.2), 
ii) {\it an enumeration formula of the chambers in the principal 
cones in terms of $\Gamma$} (Theorem 4.1) and 
iii) {\it an explanation of the formula by the block decomposition 
of the principal cone} (Cororally to Theorem 5.1).
 
When $\Gamma$ is a Coxeter diagram $\Gamma(W)$ of 
a finite Coxeter group $W$, the principal $\Gamma(W)$-cone 
was introduced in [S1] 
in the study of the real bifurcation set, where the origin and the applications 
(e.g.\ homotopy groups of complex configuration spaces,
topological types of Morsification, see \S6 and [S1,3] for more details) 
for the present study lie. 

The principal $\Gamma$-cone for an arbitrary tree 
$\Gamma$, introduced in the present article, is 
its generalization.
As we shall see, the three main results 
can be formulated and proven purely in terms of the tree $\Gamma$ 
and not of the group  $W$ 
(e.g.\ \S5 {\it Example}). 
Therefore, in the present article,
we develop the general framework for the principal $\Gamma$-cone 
and their chambers 
in a purely combinatorial manner separated from [S1].

\medskip
The contents of the present paper are as follows. 
In \S1, we fix the basic notation 
related to $\Gamma$-cones for any graph $\Gamma$.
In \S2, we prepare two assertions to count the number of 
chambers in a $\Gamma$-cone.  
After \S3, we assume that $\Gamma$ is a tree. 
In \S3, we introduce the two principal orientations on $\Gamma$ 
and the associated two principal $\Gamma$-cones $E_\Gamma$, 
and prove the first main Theorem of the present paper: 
{\it the principal $\Gamma$-cones are 
the $\Gamma$-cones which contain strictly 
maximal number of chambers}. 
In \S4, as the second main Theorem, we give the formula enumerating 
chambers in each of the principal $\Gamma$-cone in terms of $\Gamma$. 
{\it The formula is a sum of a finite number of terms where each term is a 
hook length formula for a rooted tree.}
We give an alternative proof of the formula in \S5, where 
{\it this finite sum formula 
is geometrically explained by the finite block decomposition 
of the principal $\Gamma$-cone.}
In \S6, we explain a relation of the principal $\Gamma(W)$-cone with 
the real bifurcation set of a Coxeter group $W$ [S1], which
motivated the present study. 
At the end of \S6, we compare $\Gamma$-cones with a somewhat similar concept: 
Springer cones [Ar1][Sp], and clarify the relationship between the two.

The author is grateful to Timothy Logvinenko for useful suggestions 
and to Victor Reiner for reminding some terminologies 
and references in combinatorics and, in particular, the hook length formula [K1-2].

\section{The $\Gamma$-cones and their chamber decomposition}
For a finite graph $\Gamma$, we fix notation and terminology of $\Gamma$-cones  
and their chamber decomposition.
They naturally correspond to some combinatorial 
structures  on the graph $\Gamma$ 
(e.g.\ [G-Z], [St1]). 


\medskip
Let $\Pi$ be a finite set with $\#\Pi\!=\!l\!\ge\!1$. 
A {\it configuration} in $\R$ labeled by $\Pi$ is a map $v\!:\!\Pi\!\to\! \R$ up 
to the equivalence by the translation automorphism of $\R$.
The set of all configurations labeled by $\Pi$ is given by the vector space
of rank $l\!-\!1$:
\vspace{-0.1cm}
\begin{equation}
\vspace{-0.1cm}
\label{eq:}
V_\Pi:= \oplus_{\al\in\Pi} \R v_\al/\R\cdot v_\Pi
\end{equation}
where $v_\al\ (\al\!\in\!\Pi)$ is the delta function on $\Pi$ at $\al$ 
(i.e.\ $v_\al(\beta)\!=\!\delta_{\al\beta}$ $\forall\al,\beta$) and 
$v_\Pi\!:=\!\sum_{\al\in\Pi}v_\al$ is the unit function on $\Pi$
(i.e.\ $v_\Pi(\al)\!=\!1$ $\forall\al$). 
The permutation group $\mathfrak{S}(\Pi)$ acts on $\{v_\al\}_{\al\in\Pi}$ 
fixing $v_\Pi$, and, hence, the action extends linearly to $V_\Pi$.
Let $\{\lambda_\al\}_{\al\in\Pi}$ be the dual basis of 
$\{v_\al\}_{\al\in\Pi}$, 
so that the differences 
$\lambda_{\al\beta}\!:=\!\lambda_\al\!-\!\lambda_\beta$ for 
$\al,\beta\!\in\!\Pi$ 
are well-defined linear forms on $V_\Pi$ 
(forming the root system of type $A_{l-1}$). 
The zero locus $H_{\al\beta}$ of $\lambda_{\al\beta}$ 
($\al\!\not=\!\beta$) in $V_\Pi$ is the reflection 
hyperplane of the reflection action induced by the transposition 
$(\al,\beta)$. 
The complement of the union $\cup_{\al,\beta\in\Pi,\al\not=\beta}H_{\al\beta}$ in $V_\Pi$ 
decomposes into $l!$ connected components, called 
{\it chambers} (in fact, Weyl chambers of type $A_{l-1}$,
see {\it Remark} 1.\ below). 
The set of all chambers is naturally bijective to the set 

\newpage
\noindent
$Ord(\Pi)\!:=\!\{$all {\it linear orderings} on the 
set $\Pi\}$ by the correspondence:

\vspace{0.1cm}
\centerline{
 $c:=\{\al_1\!<_c\!\ldots\!<_c\!\al_l\}\in Ord(\Pi) \ \leftrightarrow \
C_c:=\cap_{i=1}^{l-1} \{v\in V_\Pi\mid \lambda_{\al_{i}\al_{i+1}}(v)\!<\!0\}$.
}
\vspace{0.1cm}

\noindent
Here, the order-relation with respect to $c\in Ord(\Pi)$ is denoted by 
$<_c$, and the corresponding chamber is denoted by $C_c$. 
If we denote by $-c$ the reversed ordering of $c$, then one has 
$C_{-c}=-C_c$.

A {\it graph} $\Gamma$ on $\Pi$ is 
a one-dimensional simplicial complex 
whose set of vertices is $\Pi$.
An edge connecting vertices $\al$ and $\beta$ (if it exists) 
is denoted by $\overline{\al\beta}=\overline{\beta\al}$.
The set of all edges of $\Gamma$ is denoted by $Edge(\Gamma)$. 
By an abuse of notation, we shall sometimes denote the set of vertices by
$|\Gamma|$, and write ``a vertex $\al\in\Gamma$''
instead of ``a vertex $\al\in \Pi$''.

\medskip
\noindent
{\bf Definition.}  
A {\it $\Gamma$-cone} is a connected component of 
$V_\Pi\setminus \cup_{\overline{\al\beta}\in Edge(\Gamma)}H_{\al\beta}$.
By definition, each $\Gamma$-cone is subdivided into 
chambers. 

\begin{remark}
One should not confuse the chambers in the present paper, 
which are of type $A_{l-1}$ living in the $(l\!-\!1)$-dimensional space $V_\Pi$,
with the chambers of the Weyl group 
$W(\Gamma)$ of the graph $\Gamma$   ([B, VI,1.5]), which 
live in the $l$-dimensional space $\oplus \R\Pi$ (see also \S6 {\it Note} 2).\!
\end{remark}
\vspace{-0.1cm}
The $\Gamma$-cones and the chambers contained in a 
$\Gamma$-cone are described using {\it acyclic orientations} on 
$\Gamma$ defined as follows:
a collection $o$ of orientations 
$\al\!<_o\!\beta$ for all edges $\overline{\al\beta}\!\in\! Edge(\Gamma)$ 
is called {\it acyclic} if the oriented graph $(\Gamma,o)$ does not 
contain an oriented cycle (see [St2]). 
In the present paper, we shall consider only acyclic orientations. 

An acyclic orientation $o$ on $\Gamma$ defines, by transitive closure, 
a partial ordering on $\Pi$, which, for an abuse of notation, we shall 
denote by $\le_o$. 
(however, $\al\!<_o\!\beta$, i.e.\ $\al\!\le\!_o\beta$ and $\al\!\not=\!\beta$, 
may not imply the existence of an edge $\overline{\al\beta}$).
We shall denote $o_1\!\le\! o_2$ 
if $\al\!<_{o_1}\!\beta$ implies  $\al\!<_{o_2}\!\beta$.
Put
\vspace{-0.1cm}
\begin{equation}
Or(\Gamma):=\{\text{all acyclic orientations on }\Gamma \}.
\end{equation}
\vspace{-0.1cm}
The followings are immediate consequences of the definition
and are well-known [G-Z].
For a convenience of the reader, we sketch the proofs.
\begin{assertion}
{\rm 1.} For an orientation $o\in Or(\Gamma)$, define a cone:
\begin{equation}
E_o:=\cap_{ \overline{\al\beta}\in Edge(\Gamma) 
\text{ oriented as } \al <_o \beta}
\{v\in V_{\Pi}\mid \lambda_{\al\beta}(v)\!<\!0\}.
\end{equation}
Then $E_o$ is a $\Gamma$-cone.
The  correspondence $o \mapsto E_o$ induces a bijection 
\begin{equation}
Or(\Gamma)\simeq\{\Gamma\text{-cones}\}.
\end{equation}

{\rm 2.} A chamber $C_c$ for $c\in Ord(\Pi)$ is contained in 
the $\Gamma$-cone $E_o$ for  $o\in Or(\Gamma)$ 
if and only if $c$ is a linear extension of $o$, i.e.\ $o=c|_{Edge(\Gamma)}$. 

\medskip
{\rm 3.} The reflection hyperplane $H_{\al\beta}$ for $\al,\beta\in\Pi$ intersects 
with the cone $E_o$ if and only if $\al$ and $\beta$ are disconnected 
in $\Gamma$ and the orientation $o$ on $\Gamma$ induces an acyclic 
orientation
on the quotient graph $\Gamma/\!\sim$, where 
$\Gamma/\!\sim$ is obtained from $\Gamma$ by the identification of $\al$ and $\beta$.
\end{assertion}
\begin{proof}
1. For any $o\in Or(\Gamma)$, let us show that $E_o\not=\emptyset$,
that is: there exists a map 
$v:\Pi\!\to\! \R$ such that $v(\al)\!<\!v(\beta)$ if $\al\!<_o\!\beta$. 
We prove this by an induction on $l$.
Acyclicity of $(\Gamma,o)$ implies an existence of a 
{\it minimal} vertex 
$\al\!\in\! \Pi$: for any edge $\overline{\al\beta}\!\in\! Edge(\Gamma)$, 
one has $\al\!<_o\!\beta$. Put $\Pi'\!:=\!\Pi\setminus \{\al\}$. 
Then clearly $o'\!:=\!o|_{\Pi'}$ is an orientation on the graph 
$\Gamma'\!:=\!\Gamma|_{\Pi'}$. By the induction hypothesis, 
there exists a map $v'\!:\!\Pi'\!\to\!\R$ preserving the sub-orientation $o'$.
Then, $v$ is an extension of $v'$ by choosing the value $v(\al)$
from the non-empty set 
$\R\setminus\! \cup_{\beta\in\Pi', \overline{\al\beta}\in Edge(\Gamma)}
[v'(\beta),\infty)$.

Conversely, for a given $\Gamma$-cone $E$, define the 
edge $\overline{\al\beta}\in Edge(\Gamma)$ to have the orientation 
$\al\!<_E\!\beta$  if 
$\lambda_{\al\beta}|_{E}<0$. This defines an acyclic orientation $o_E$ on $\Gamma$.
These establish the bijection (4).

2. The inclusion $C_c\!\subset\! E_o$ is equivalent 
to the inclusions 
$C_c\!\subset\!\{\lambda_{\al\beta}\!<\!0\} \Leftrightarrow 
\al\!<_c\!\beta$ for any oriented edge 
$\overline{\al\beta}$ with $\al\!<_o\!\beta$. 

3. That $H_{\al\beta}\!\cap\! E_o\!\not=\!\emptyset$ is equivalent that there exists a
map $v\!:\!\Pi\!\to\!\!\R$ satisfying $v(\al)\!=\!v(\beta)$ and the inequalities (3). 
Apply again the argument in 1. to this setting.
\end{proof}

According to the previous assertion, we put
 \begin{equation}
\Sigma(o):=\{\ c\in Ord(\Pi)\ \mid\ o=c|_{Edge(\Gamma)}\ \},\\
 \end{equation}
and identify $\Sigma(o)$ with the set of chambers contained in $E_o$. 
Let us introduce a numerical invariant for the orientation $o\in Or(\Gamma)$:
\begin{equation}
\label{eq:}
\quad \sigma(E_o):=\sigma(o):=\# \Sigma(o)=\#\{\text{chambers contained in } E_o\}.\!\!\!\!\!\!
\end{equation}

If we denote by $-o$ the reversed orientation of $o$, one has 
$E_{-o}=-E_o$ and, therefore, $\Sigma(-o)=-\Sigma(o)$ and 
$\sigma(-o)=\sigma(o)$.


\medskip
In order to obtain the smallest oriented graph, which gives the same cone, 
we introduce the 
{\it reduced oriented graph} $o_{red}$, 
or the {\it Hasse diagram}: 
an oriented edge $\al\!<_o\!\beta$ of the oriented graph $o$ is called 
{\it reducible}
if there is a sequence 
$\al_0\!=\!\al,\al_1,\cdots,\al_{k-1},\al_k\!=\!\beta\in\Pi$ for some  
$k\!\in\!\Z_{>1}$ such that 
$\al_{i-\!1}\!\!<_o\!\!\al_i$ for $i\!=\!1,\!\cdots\!,k$.
Then the {\it reduction $o_{red}$ of $o$}, also called 
the {\it Hasse diagram} of $o$, is uniquely defined from 
$o$ by

\quad $o_{red}:=$ 
the oriented graph obtained from $o$ by deleting 

\qquad \qquad\ all reducible edges.

One easily observes that 
i) {\it the associated cones coincide: $E_o\!=\!E_{o_{red}}$,}
ii) {\it there is a bijection between the edges of 
$o_{red}$
and the $(l\!\!-\!2)$-dimensional faces of $E_o\!=\!\!E_{o_{red}}$.}\
Consequently,\ one has: 
iii) {\it  $E_o\!=\!E_{o'}$ for oriented graphs $o$ and $o'$ on $\Pi$, 
if and only if $o_{red}\!=\!o'_{red}$.}

\begin{remark} As described above, the geometry of 
the $\Gamma$-cones in $V_{\Pi}$ has a natural relation  with 
the combinatorics (partially ordered structures) on the set $\Pi$. 
The enumeration of $\sigma(o)$ 
is {\it the basic problem of enumeration of linear extension of a partially ordered 
set} (e.g., see [St1]).

On the other hand, a $\Gamma$-cone with the 
subdivision into chambers 
appears naturally in the study of finite Coxeter group $W$ as follows, 
where  
$\Pi$ stands for a {\it simple generator system} of $W$. Then,  
a linear ordering $\al_1\!\!<_c\!\!\al_2\!\!<_c\!\!\cdots\!\!<_c\!\!\al_l$ 
of $\Pi$ defines 
a {\it Coxeter element} $\al_1\cdots\al_l\!\in\! W$. Two Coxeter elements
coincide if the corresponding chambers belong to the same $\Gamma(W)$-cone 
for the Coxeter-Dynkin diagram $\Gamma(W)$ on $\Pi$.
The principal $\Gamma(W)$-cone $E_{\Gamma(W)}$,  
which we shall introduce in \S3,  
has the particular geometric significance, 
for which we refer to [S1] (see \S7). 

It may be worthwhile to mention 
that a choice of the orientations on the Coxeter-Dynkin diagram
$\Gamma(W)$ plays often an important role in the studies related to 
root systems
(recent examples: a stability condition on a triangulated category [K-S-T], 
the quantized Toda equations [E]).
\end{remark}



\vspace{-0.3cm}

\section{A decomposition formula}
We prepare two Assertions which are to help calculating $\sigma(o)$.
They are used in the proof of the Theorems in \S3,4 and 5. 
The reader may skip this section at first reading.

\medskip
For any $o\in Or(\Gamma)$, $\al\in\Pi$ and $r\in\Z_{\ge0}$, 
we put 
\begin{align}
\Sigma(o,\al,r)&:=\{c\in \Sigma(o)\mid \#\{\beta\in\Pi\mid \al<_c\beta\}=r\},\\
\sigma(o,\al,r)&:=\#\Sigma(o,\al,r).
\end{align}
Obviously, one has the disjoint decomposition  
$\Sigma(o)=\coprod_{r=0}^{l-1}\Sigma(o,\al,r)$ for any $\al\in\Pi$ 
so that
$\sigma(o)=\sum_{r=0}^{l-1}\sigma(o,\al,r)$.

1. 
Let $\Gamma_1,\cdots,\Gamma_k$ ($k\in\Z_{>0}$) be graphs, which contain 
the unique common vertex $\al$. Let us denote by 
\begin{equation}
\begin{array}{ll}
\Gamma_1 \coprod_\al\cdots\coprod_\al \Gamma_k,
\end{array}
\end{equation}
the graph obtained by the disjoint union of 
the graphs $\Gamma_i$ ($i\!=\!1,\cdots,k$)
up to an identification of the common vertex $\al$. 

\begin{assertion} Let 
$\Gamma$ be a graph decomposing as (9).
For an orienta-

\noindent
tion  $o\!\in\! Or(\Gamma)$, 
put $o_i\!:=\!o|_{\Gamma_i}\!\in\! Or(\Gamma_i)$ $(1\!\!\le\!\! i\!\!\le\!\!k)$ 
and $l_i\!:=\!\#\Gamma_i$ $(1\!\!\le\!\! i\!\!\le\!\!k)$ 
so that $l\!=\!\sum_il_i\!-\!k\!+\!1$. 
Then, for any $r\!\in\!\Z_{\ge0}$, one has the following:
\begin{equation}
\sigma(o,\al,r)=\!\!\! \sum_{\substack{r_1,\cdots,r_k\in\Z_{\ge0}\\
r_1+\cdots+r_k=r}}\!\!
\sigma(o_1,\al,r_1)\cdots\sigma(o_k,\al,r_k)(_{\ r_1,\cdots,r_k}^{r_1+\cdots+r_k})
(_{l_1-r_1-1,\cdots,l_k-r_k-1}^{\ \ \ l-1-r_1-\cdots-r_k}),
\end{equation}
where  
$(_{\ r_1,\cdots,r_k}^{r_1+\cdots+r_k})\!:=\!(r_1\!+\!\cdots\!+\!r_k)!/r_1!\cdots r_k!$ is the 
multinomial coefficient. 
By summing the formula {\rm (10)} for all $r\in\Z_{\ge0}$, one obtains the following:
\begin{equation}
\sigma(o)=\!\!\! \sum_{r_1,\cdots,r_k\in\Z_{\ge0}}\!\!
\sigma(o_1,\al,r_1)\cdots\sigma(o_k,\al,r_k)(_{\ r_1,\cdots,r_k}^{r_1+\cdots+r_k})
(_{l_1-r_1-1,\cdots,l_k-r_k-1}^{\ \ \ l-1-r_1-\cdots-r_k}).
\end{equation}
\end{assertion}
\begin{proof}
 It is sufficient to prove only (10).

Consider the projection 
$\Sigma(o)\to \Sigma(o_1)\times \cdots \times \Sigma(o_k),$ 
$c\mapsto (c|_{\Gamma_i})_{i=1,\cdots,k}$.
The projection decomposes into projections 
\begin{equation}
\Sigma(o,\al,r)\to \coprod_{\substack{r_1,\cdots,r_k\in\Z_{\ge0}\\
r_1+\cdots+r_k=r}}\Sigma(o_1,\al,r_1)\times \cdots \times \Sigma(o_k,\al,r_k)
\end{equation}
for $r\in\Z_{\ge0}$.
Let us see that the cardinality of the inverse image of a point 
$(c_1,\cdots,c_k)\in \Sigma(o_1,\al,r_1)\times \cdots \times \Sigma(o_k,\al,r_k)$ 
is equal to the product of two combination number: 
$(_{\ r_1,\cdots,r_k}^{r_1+\cdots+r_k})(_{l_1-r_1-1,\cdots,l_k-r_k-1}^{\ \ \ l-1-r_1-\cdots-r_k})$.

Put $\Gamma_i^+\!:=\!\{\beta\!\in\!|\Gamma_i| \mid \al\!<_{c_i}\!\beta\}$ and 
$\Gamma_i^-\!:=\!\{\beta\!\in\!|\Gamma_i| \mid \beta\!<_{c_i}\!\al\}$.
Then, an ordering $c\in\Sigma(o,\al,r)$ 
is in the inverse image, if the $r\!=\!r_1\!+\!\cdots\!+\!r_k$ elements 
$\amalg_{i=1}^k\Gamma_i^+$ lie in the right hand side of $\al$ 
and the
$l\!-\!r\!-\!1\!=\!(l_1\!-\!r_1\!-1)\!+\!\cdots\!+\!(l_k\!-\!r_k\!-\!1)$ 
elements $\amalg_{i=1}^k\Gamma_i^-$
lie in left hand side of $\al$ with respect to $c$,  
and the restrictions  $c|_{\Gamma_i^\pm}$ 
are equal to the pre-fixed linear orderings $c_i^\pm:=c_i|_{\Gamma_i^\pm}$ for $i=1,\cdots,k$.
Thus the set of $c$ in the inverse image is bijective to  
$\Sigma(\amalg_{i=1}^kc_i^+) \times \Sigma(\amalg_{i=1}^kc_i^-)$, 
where $\amalg_{i=1}^kc_i^\pm$ is the partial ordering structure on 
$\Gamma^\pm:=\amalg_{i=1}^k\Gamma_i^\pm$. 
Since $\Sigma(\amalg_{i=1}^kc_i^+)$ is just the set of shuffles of 
$k$ sets of cardinalities $r_1$,\ldots,$r_k$, its cardinality is given by the 
combination number, i.e.\  
$\sigma(\amalg_{i=1}^kc_i^+)\!=(_{\ r_1,\cdots,r_k}^{r_1+\cdots+r_k})$.

Similarly, one has
$\sigma(\amalg_{i=1}^kc_i^-)\!=\!(_{l_1-r_1-1,\cdots,l_k-r_k-1}^{\ \ \ l-1-r_1-\cdots-r_k})
$.
\end{proof}

2. \ A vertex $\al\in \Pi$ 
is called {\it maximal} (resp.\ {\it minimal}) 
with respect to $o\in Or(\Gamma)$, 
if  $\beta\!<_o\!\al$ (resp.\ $\al\!<_o\!\beta$) 
for any edge $\overline{\al\beta}\in Edge(\Gamma)$ at $\al$.

\begin{assertion} 
\noindent
If $\al$ is maximal with respect to $o$, then one has 

\centerline{
$\sigma(o,\al,0)\ge\sigma(o,\al,1)\ge\cdots\ge 
\sigma(o,\al,l\!-\!2)\ge\sigma(o,\al,l\!-\!1)$.
}

\noindent
If $\al$ is minimal with respect to $o$, then one has 

\centerline{
$\sigma(o,\al,0)\le\sigma(o,\al,1)\le\cdots\le 
\sigma(o,\al,l\!-\!2)\le\sigma(o,\al,l\!-\!1)$.
}
\noindent
If $\al$ is non-isolated in $\Gamma$, the smallest 
terms in the sequences are zero.
\end{assertion}
\begin{proof} 
We show only the first case. The latter case is shown similarly.

It is sufficient to show that there is an injection map 
$\Sigma(o,\al,r)\to\Sigma(o,\al,r\!-\!1)$ for $r>0$.
In fact the map is constructed as follows:
let $c=\{A<_c\al<_c\beta<_c B\}\in \Sigma(o,\al,r)$ where $\beta\in \Pi$ 
and $A$ and $B$ 
are such linear sequence of inequalities of elements of $\Pi$ that 
the length of $B$ is equal to $r-1$ (this is possible 
since $r\ge1$). Then we set 
$c':=\{A<_c\beta<_c\al<_c B\}\in \Sigma(o,\al,r-1)$ where $c'$ 
is well defined since $\al$ is maximal. The correspondence 
$c\mapsto c'$ is clearly injective.  

If $\al$ is non-isolated, then the set $\Sigma(o,\al,l\!-\!1)$ is empty, since 
there exists a vertex $\beta\in\Pi$ such that 
$\overline{\beta\al}\in Edge(\Gamma)$ and $\beta<_o\al$ and 
hence for any $c\in\Sigma(o)$ one has $\beta<_c\al$ and $c\not\in \Sigma(o,\al,l-1)$.
\end{proof}

\section{Principal $\Gamma$-cones}
A graph $\Gamma$ is called a {\it tree}
if it is connected and has no cycle. 
For a tree $\Gamma$, we introduce two particular $\Gamma$-cones, 
called the principal $\Gamma$-cones.
The first main result of the present paper, 
formulated in Theorem (3.2), is to 
characterize\! the\! principal\! 
$\Gamma$-cones.

\medskip
The following is a characterization of trees in terms of $\Gamma$-cones.

\begin{assertion}
 Let $\Gamma$ be a graph on $\Pi$. 
Then, 
$\{H_{\al\beta}\}_{\overline{\al\beta}\in Edge(\Gamma)}$ 
forms a system of coordinate hyperplanes of $V_\Pi$ 
if and only if $\Gamma$ is a tree. 
\end{assertion}

 \begin{proof}
 For each edge $\overline{\al\beta}$ of $\Gamma$, we choose one of 
$\lambda_{\al\beta}$ or $\lambda_{\beta\al}$.
Then, it is immediate that 
 i) $\{\lambda_{\al\beta}\}_{\overline{\al\beta}\in Edge(\Gamma)}$ 
 is linearly independent if and only if  $\Gamma$ 
 does not contain a cycle, and   
 ii)  $\{\lambda_{\al\beta}\}_{\overline{\al\beta}\in Edge(\Gamma)}$
 spans the dual space of
 $V_\Pi$ if and only if  $\Gamma$ is connected.
 \end{proof}

Consequently, a $\Gamma$-cone $E_o$ is simplicial (i.e.\ the cone over a simplex) 
if and only if $o_{red}$ is a tree.
From now on in the present paper, we shall assume that $\Gamma$ 
is a tree on $\Pi$.  Then the system of coordinate hyperplanes 
 $\{H_{\al\beta}\}_{\overline{\al\beta}\in Edge(\Gamma)}$ 
cuts the vector space 
 $V_\Pi$ into $2^{l-1}$-number of orthants, 
each of which is simplicial. 
A question is a characterization of
the $\Gamma$-vector: $(\sigma(o))_{o\in Or(\Gamma)}$ 
of size $2^{l-1}$. 
One distinguished property of
the $\Gamma$-vector is that it contains a unique 
(up to change the sign of orientations, see {\it Note.}\  below) 
maximal entry, which we explain now.

\begin{definition} Let $\Gamma$ be a tree (or more generally, a connected graph).

1. A {\it principal decomposition of $\Gamma$} is  
the ordered pair $\{\Pi_1,\Pi_2\}$ of subsets of $\Pi$ such that 
i) one has the disjoint decomposition:
\begin{equation}
\Pi=\Pi_1\amalg \Pi_2,
\end{equation}
ii) each $\Pi_i$ is totally disconnected (discrete) in $\Gamma$.  

2. A {\it principal orientation}
 on $\Gamma$ is an element in $Or(\Gamma)$ attached to 
a principal decomposition $\{\Pi_1,\Pi_2\}$ as follows:
\begin{equation}
o_{\Pi_1,\Pi_2}\!:=\!\{\ \al<_{o_{\Pi_1,\Pi_2}}\beta 
\ \text{ for } \overline{\al\beta}\!\in\! Edge(\Gamma)  \text{ with }
 \al\!\in\!\Pi_1, \beta\!\in\!\Pi_2 \}.\!\!\!\!\!\!\!
\end{equation} 

3.  A {\it principal $\Gamma$-cone}  
is the $\Gamma$-cone attached to a principal orientation 
$o_{\Pi_1,\Pi_2}$. That is: 
\begin{equation}
\begin{array}{lll}
&&E_{\Pi_1,\Pi_2}\ :=\ E_{o_{\Pi_1,\Pi_2}} \\
&=\!\!&\{v\in V_\Pi\mid 
\lambda_{\al\beta}(v)\!<0\! \text{ for } \overline{\al\beta}\!\in\! Edge(\Gamma)  \text{ with }
 \al\!\in\!\Pi_1, \beta\!\in\!\Pi_2\}.\!\!\!\!\!
\end{array}
\end{equation}
\end{definition} 

It is easy to see that there exist exactly  
two principal decompositions for any tree $\Gamma\!\not=\!\emptyset$,  
and that if one is given by $\{\Pi_1,\Pi_2\}$ 
then the other one is given by $\{\Pi_2,\Pi_1\}$. 
For simplicity, we shall confuse the expression (13) with 
the principal decomposition $\{\Pi_1,\Pi_2\}$.

4. Since 
 $o_{\Pi_2,\Pi_1}\!=\!-o_{\Pi_1,\Pi_2}$ and
$E_{\Pi_2,\Pi_1}\!=\!-E_{\Pi_1,\Pi_2}$, two principal 
$\Gamma$-cones (see {\it Note.} below)
are isomorphic to each other as abstract cones by the multiplication of $-1$. 
The isomorphism preserves the subdivisions into chambers.
The isomorphism class of 
the pair of the cones and its subdivision into chambers 
is called the {\it  principal $\Gamma$-cone}.
For an abuse of notation, it is denoted by
\begin{equation}
E_\Gamma:= E_{\Pi_1,\Pi_2}\simeq E_{\Pi_2,\Pi_1}.
\end{equation} 

\begin{note} 
If $\Gamma$ is $\Gamma(A_1)\!=$ one point graph, then $\Gamma$ has only 
one orientation, say $o_{A_1}$. 
Therefore, the two principal orientations 
$o_{\Pi_1,\Pi_2}$ and $o_{\Pi_2,\Pi_1}$ coincide with 
$o_{A_1}$, i.e.\ $o_{\Pi_1,\Pi_2}\!=\!o_{\Pi_2,\Pi_1}\!=\o_{A_1}$. 
The $V_\Pi$ is the zero vector space $\{0\}$, which is equal to 
the only principal cone $E_{\Gamma(A_1)}\!=\!\{0\}$. The cone consists of 
single chamber $C\!:=\!\{0\}$, i.e.\  
$\Sigma(o_{A_1})\!=\!\{C\}$ and $\sigma(o_{A_1})\!=\!1$. 
{\it Except for this case, there are always 
two principal $\Gamma$-cones, which 
are in the opposite position with respect to the origin in $V_\Pi$.}
\end{note}

\medskip
The following is the first main theorem of the present paper,
which characterizes the principal $\Gamma$-cone and justifies 
its naming ``principal''.

\begin{theorem}
Let $\Gamma$ be a tree on $\Pi$. 
The principal $\Gamma$-cone $E_\Gamma$ is the $\Gamma$-cone which contains
the strictly maximal number of chambers.
That is: a $\Gamma$-cone $E_o$ for $o\!\in\! Or(\Gamma)$ 
is principal, i.e.\ $o$ is one of the two orientations 
$o_{\Pi_1,\Pi_2}$ or $o_{\Pi_2,\Pi_1}$,  
if and only if $\sigma(o)\!=\!\sigma(\Gamma)$, where 
\begin{equation}
\sigma(\Gamma):=\max\{\sigma(p)\mid p\in Or(\Gamma)\}.
\end{equation}
\end{theorem}
\begin{proof}
With the results established in \S2, the proof is straightforward:
suppose $o\in Or(\Gamma)$ is not principal, that is: 
there exist $\al,\beta,\gamma\in\Pi$ 
with $\gamma<_o\al<_o\beta$. 
Actually, $\Gamma$ decomposes as $\Gamma=\Gamma_+\coprod_\al\Gamma_-$, 
where $\Gamma_+$ (resp.\ $\Gamma_-$) is a full subgraph of $\Gamma$ 
containing $\al$ and any connected component of $\Gamma\setminus\{\al\}$
which contains a vertex $\beta$ s.t. $\al<_o\beta$ (resp.\ $\al>_o\beta$).
By the assumption on $o$, 
one has $\Gamma_\pm\not=\emptyset$.

Put $o_+:=o|_{\Gamma_+} \in Or(\Gamma_+)$ and $o_-:=o|_{\Gamma_-} \in Or(\Gamma_-)$. 

\begin{assertion} Define a new orientation 
$\tilde o\!\in\! Or(\Gamma)$ by the following rule: 
$\tilde o$ agrees with $o_+$ on  $\Gamma_+$ 
and with $-o_-$ on $\Gamma_-$. 
Then $\sigma(\tilde o)>\sigma(o)$.
\end{assertion}
\begin{proof}
For a proof of the Assertion, we apply the formula (11) in Assertion 2.1
to the decomposition $\Gamma=\Gamma_+\coprod_\al\Gamma_-$and  
to $o,\tilde o\in Or(\Gamma)$:
\begin{align*}
\sigma(o) = \sum_{r_+=0}^{l_+}\sum_{r_-=0}^{l_-}
\sigma(o_+,\al,r_+)\sigma(o_-,\al,r_-) C_{r_+,\ r_-}C_{l_+-r_+,\ l_--r_-},\\
\sigma(\tilde o) = \sum_{r_+=0}^{l_+}\sum_{r_-=0}^{l_-}
\sigma(o_+,\al,r_+)\sigma(o_-,\al,r_-) C_{r_+,\ l_--r_-}C_{l_+-r_+,\ r_-},
\end{align*}
where $l_+:=\#\Gamma_+\!-\!1>0$,\ $l_-:=\#\Gamma_-\!-\!1>0$ 
such that $l=l_++l_-+1$,  
and the notation $C_{r_+,\ r_-}$ 
means the binomial coefficients 
$\frac{(r_++r_-)!}{r_+!\ r_-!}$.

We want to calculate the difference $\sigma(\tilde o)-\sigma(o)$ term-by-term. 
Observe that the terms for $r_+=l_+/2$ (if $l_+$ is even) and
the terms for $r_-=l_-/2$ (if $l_-$ is even) in the two formulae  give 
the same value and so cancel each other in the difference. 
Therefore, we decompose the region $[0,l_+]\times [0,l_-]$ 
of the summation index $(r_+,r_-)$ into 4 regions 
according to whether $r_+$ is larger or less than $l_+/2$ and 
$r_-$ is larger or less than $l_-/2$.

For an  index $(r_+,r_-)$ in the region $[0,l_+/2)\times [0,l_-/2)$, 
we consider 4 indices $(r_+,r_-)$,  $(r_+,r_-^*)$,  $(r_+^*,r_-)$ 
and  $(r_+^*,r_-^*)$ in the 4 regions simultaneously, where 
$r_+^*:=l_+-r_+$ and $r_-^*:=l_--r_-$. 
Let us explicitly write down the difference between these 4 terms in 
$\sigma(\tilde o)$ and in $\sigma(o)$:
\[\begin{array}{lll}
&\! \sigma(r_+)\sigma(r_-) C_{r_+,r_-^*}C_{r_+^*,r_-}
+\sigma(r_+)\sigma(r_-^*) C_{r_+,r_-}C_{r_+^*,r_-^*}\\
+\!&\!\!\sigma(r_+^*)\sigma(r_-) C_{r_+^*,r_-^*}C_{r_+,r_-}
+\sigma(r_+^*)\sigma(r_-^*) C_{r_+^*,r_-}C_{r_+,r_-^*}\\
-\!&\!\!\sigma(r_+)\sigma(r_-) C_{r_+,r_-}C_{r_+^*,r_-^*}
-\sigma(r_+)\sigma(r_-^*) C_{r_+,r_-^*}C_{r_+^*,r_-}\\
-\!&\!\!\sigma(r_+^*)\sigma(r_-) C_{r_+^*,r_-}C_{r_+,r_-^*}
-\sigma(r_+^*)\sigma(r_-^*) C_{r_+^*,r_-^*}C_{r_+,r_-},
\end{array}
\]
where we used the simplified notation 
$\sigma(r_+):=\sigma(o_+,\al,r_+)$, $\sigma(r_+^*):=\sigma(o_+,\al,r_+^*)$,
$\sigma(r_-):=\sigma(o_-,\al,r_-)$ and $\sigma(r_-^*):=\sigma(o_-,\al,r_-^*)$.

One can factorize this difference as follows:
\[
(\sigma(r_+)-\sigma(r_+^*))(\sigma(r_-^*)-\sigma(r_-))
( C_{r_+,r_-}C_{r_+^*,r_-^*}- C_{r_+,r_-^*}C_{r_+^*,r_-}).
\]

Let us examine the sign of the factors and demonstrate that 
the product turns out to 
be non-negative.
First, recall that the vertex $\al$ is minimal in $\Gamma_+$ and maximal in 
$\Gamma_-$ by definition. Note also that 
$r_+<l_+/2<r_+^*$ and  $r_-<l_-/2<r_-^*$.
Therefore, applying \S2 Assertion 2.2, we observe that 
$(\sigma(r_+^*)-\sigma(r_+))(\sigma(r_-)-\sigma(r_-^*))\ge0$. 
Next, let us examine the last factor. For this purpose, 
we use the proportion of the two terms in the last factor:
\[
\frac{C_{r_+,r_-}C_{r_+^*,r_-^*}}{C_{r_+,r_-^*}C_{r_+^*,r_-}}
=\frac{(r_+^*+r_-^*)!}{(r_++r_-^*)!}\cdot \frac{(r_++r_-)!}{(r_+^*+r_-)!}.
\]
Using the fact that $r_+<r_+^*$, 
one has  $r_+^*+r_-^*>r_++r_-^*$ and $r_++r_-<r_+^*+r_-$. 
Hence, the expression can be reduced to 
$
\prod_{k=r_++1}^{r_+^*} \frac{r_-^*+k}{r_-+k}
$,
where each factor is larger than 1 since $r_-+k<r_-^*+k$ and the number 
of the factors is $r_+^*-r_+>0$ so that the result is always larger than $1$. 
These together imply that the difference of the 4 terms is non-negative.

By summing up terms for all indices $(r_+,r_-)$ in the region 
$[0,l_+/2)\times [0,l_-/2)$, we see that the difference 
$\sigma(\tilde o)-\sigma(o)$ is non-negative. To 
show that it is strictly positive, let us calculate the term for 
$(r_+,r_-)=(0,0)$.
Then, \S2, Assertion 2.2 again, one has 
$\sigma(r_+)\!=\!\sigma(r_-^*)\!=\!0$. 
Since $l_+,l_->0$ (non-principality of $\sigma$), one obtains a rather big number:

\centerline{
$\sigma(o_+,\al,l_+)\sigma(o_-,\al,0)(C_{l_+,l_-}\!-1)\!\not=\!0$.}

This completes the proof of the Assertion.
\end{proof}
The Assertion says that if an orientation $o$ on $\Gamma$ is not principal, 
it can not attain the maximal value $\sigma(\Gamma)$ of 
$\sigma(o)$ for $o\!\in\! Or(\Gamma)$. In fact, 
starting from any orientation $o\!\in\! Or(\Gamma)$, and by a 
successive application of the construction in the Assertion, 
one arrives at one of the principal orientations.
Since $E_{\Pi_1,\Pi_2}\simeq E_{\Pi_1,\Pi_2}$, one has 
$\sigma(o_{\Pi_1,\Pi_2})=\sigma(o_{\Pi_1,\Pi_2})$. This number gives 
the maximal value $\sigma(\Gamma)$. 

This completes the proof of the Theorem.
\end{proof} 

We shall call $\sigma(\Gamma)$ the {\it principal number} of $\Gamma$.
\begin{remark}
Some particular cases of Theorem 3.2 were known already.

If $\Gamma$ is a linear graph of type $A_l$, then the (principal) 
$\Gamma$-cones coincide with the (principal) Springer cone of type 
$A_{l-1}$ (see [Ar],[Sp] and the latter half of \S7 of the present paper).
In that case, the result is shown in [Sp, Prop.3],
[S-Y-Z, Theo.1.2., (2)], [N],[Br].

\end{remark}

\begin{remark}
Let $\Gamma\!=\!\amalg_{i=1}^k \Gamma_i$ be the decomposition of a forest 
into trees. For an orientation $o$ on $\Gamma$, 
put $o_i\!:=\!o|_{\Gamma_i}$. 
Since $\sigma(o)\!=\!\Pi_i\sigma(o_i) 
(^{\ \ \sum\#\Gamma_i}_{\#\Gamma_1,\cdots,\#\Gamma_k})$,
the maximal number of chambers in a $\Gamma$-cone is 
attained by the orientations $o$ such that each $o_i$ is a principal 
orientation on $\Gamma_i$.
\end{remark}

\medskip
\noindent
{\it Example.}
The principal cone $E_{\Gamma(D_4)}$ consists of 
6 chambers forming a hexagon.
The $\Gamma(D_4)$-vector is 
%
$(\sigma(o))_{o\in Or(\Gamma(D_4))}=2(6,2,2,2)$.

The principal cone $E_{\Gamma(A_4)}$ consists of 
5 chambers forming a spoon graph.
The $\Gamma(A_4)$-vector is 
$(\sigma(o))_{o\in Or(\Gamma(A_4))}=2(5,3,3,1)$.

Let $\Gamma$ be a cyclic graph of 4 vertices. Even though $\Gamma$ is not a tree,
the $\Gamma$-vector contains the maximal entry:  
$(\sigma(o))_{o\in Or(\Gamma)}=2(4,2,2,1,1,1,1)$, where $(*)$ {\it the maximal 
number is 
attained by the $\Gamma$-cones corresponding to
principal decompositions} (13). 
Actually, a principal decomposition for an arbitrary $\Gamma$ 
may not exist. However, conjecturally, $(*)$ holds 
for any connected graph $\Gamma$ which admits a principal decomposition.

\begin{remark}([S2])
Assume that a connected graph $\Gamma$ admits a principal decomposition. 
Consider the lattice $L_\Gamma$ spanned by $\Pi$ with the symmetric 
bilinear form as in the usual convention in the theory of root systems.
Then the ``Coxeter element''defined as the product of reflections attached 
to the vertex in the order of a principal order (recall {\it Remark} 1) 
is i) semi-simple of finite order, or ii) quasi-unipotent 
if and only if i) $\Gamma$ is one of the Coxeter-Dynkin diagram for 
a finite Coxeter group or ii) $\Gamma$ is one of the affine 
Coxeter-Dynkin diagram, respectively.
\end{remark}

\section{Enumeration of chambers in the principal $\Gamma$ cone}


As the second main result of the present paper, we give 
an enumeration formula (19) for the principal number $\sigma(\Gamma)$.
The formula depends on a choice of a principal decomposition (13). 
It is formulated as i) a sum  
whose summation index runs over certain equivalence classes 
$\widetilde{Ord}(\Pi_1)$ of all linear orderings on $\Pi_1$ and 
ii) each summand is the quotient of $(\#\Gamma)!$ 
by a product of cardinalities of certain subgraphs. 

We will present two proofs of the formula (19). 
In the present section, we give a direct proof  
based on the principal $\Gamma$-ordering on $\Pi$. The proof in 
\S5 is based on a decomposition of the principal $\Gamma$-cone into blocks 
attached to newly introduced rooted tree structures on $\Pi$ 
and applying the hook length formula of Knuth [K2,p70].

We start with the definition of the equivalence $\sim$ on the set $Ord(\Pi_1)$.

\medskip
Let  $d\in Ord(\Pi_1)$ be an ordering on $\Pi_1$. For $v\in\Pi_1$, put 
\begin{equation}
\begin{array}{lll}
\Gamma_{d,v}&:=&\text{the connected component of }
\Gamma\!\setminus\!\{w\!\in\!\Pi_1\mid w\!<_d\!v\} \\
&&\text{ containing }v.
\end{array}
\end{equation}
In particular, one has $\Gamma_{d,v}=\Gamma$ for the smallest element $v$ of $\Pi_1$.

\begin{definition}
Two orderings $d,d'\in Ord(\Pi_1)$ are called {\it equivalent} if  
$\Gamma_{d,v}=\Gamma_{d',v}$ for all $v\in \Pi_1$.
The equivalence class of $d$ is denoted by $\tilde d$ and  
the set of all equivalence classes is denoted by 
$\widetilde{Ord}(\Pi_1)$. 
\end{definition}


\begin{theorem} Let $\Gamma$ be a tree on $\Pi$. 
 Choose a principal decomposition {\rm (13)}. 
Then the principal number $\sigma(\Gamma)$ defined in {\rm (17)} is given by 
\begin{equation}
\sigma(\Gamma)=(\#\Gamma)! \sum_{\tilde d
\in \widetilde{Ord}(\Pi_1)}\frac{1}{\prod_{v\in \Pi_1} 
\#\Gamma_{d,v}},
\end{equation}
where the terms in RHS are well-defined since 
$\Gamma_{d,v}$ depends only on the 
equivalence class $\tilde d$ of $d\in Ord(\Pi_1)$ and on $v\in \Pi_1$.
\end{theorem}
\begin{proof} 
Before we start with the proof of the formula, we reformulate the equivalence 
$\sim$ in terms of the partial orderings on the set $\Pi_1$.

\begin{assertion}
Let a linear ordering  $d\in Ord(\Pi_1)$ on $\Pi_1$ be given.
 For any two vertices $v,v'\in \Pi_1$, one has the following:

 i) There are only three cases:
\[
\Gamma_{d,v}\cap\Gamma_{d,v'}=
\begin{cases}
\emptyset \\
\Gamma_{d,v}\\
\Gamma_{d,v'}.
\end{cases}
\]

 ii) If $\Gamma_{d,v}\cap\Gamma_{d,v'}=\Gamma_{d,v}$, 
then one has $v'\le_d v$.  

\medskip
iii) The next three conditions are equivalent:

\centerline{
a) $\Gamma_{d,v}\cap\Gamma_{d,v'}=\Gamma_{d,v}$, \quad 
b) $\Gamma_{d,v}\subset\Gamma_{d,v'}$,\quad 
c) $v\in \Gamma_{d,v'}$.
}
\end{assertion}
\begin{proof} 
i) Since $d$ is a linear ordering, we may assume $v'<_dv$. The fact 
that
 $\Gamma\!\setminus\!\{w\in \Pi_1\mid w<_d v\}\subset \Gamma\!\setminus\!\{w\in \Pi_1\mid w<_d v'\}$ 
implies that the component $\Gamma_{d,v}$ is either contained in 
the component $\Gamma_{d,v'}$ 
or they are disjoint. Accordingly, the intersection is either 
$\Gamma_{d,v}$ or an empty set.

ii)  Suppose the contrary $v'\not\le_d v$. 
Then the fact that $d$ is totally ordered implies 
$v'>_d v$. Then, by the construction, $\Gamma_{d,v'}$ cannot contain $v$.
This contradicts  the assumption 
  $\Gamma_{d,v}\cap\Gamma_{d,v'}=\Gamma_{d,v}$.

iii) The implications: a) $\Rightarrow$ b) $\Rightarrow$ c) are trivial. 
Assume c). This implies 
$\Gamma_{d,v}\cap\Gamma_{d,v'}\not=\emptyset$. Suppose, further, 
$\Gamma_{d,v}\cap\Gamma_{d,v'}\not=\Gamma_{d,v}$. Then i) implies 
$\Gamma_{d,v}\cap\Gamma_{d,v'}=\Gamma_{d,v'}\not=\Gamma_{d,v}$, 
and, hence, $\Gamma_{d,v'}\subset\not=\Gamma_{d,v}$.
This means that $\Gamma_{d,v'}$ is a connected component 
obtained by deleting 
strictly more vertices than  those for $\Gamma_{d,v}$. This is possible only 
 when $v<_d v'$. Then, $v\not\in \Gamma_{d,v'}$. 
A contradiction to the assumption c) !
\end{proof}

\begin{definition}
 To the equivalence class $\tilde d$ in $\widetilde{Ord}(\Pi_1)$  of 
$d\in Ord(\Pi_1)$, we attach a {\it partial ordering}
on $\Pi_1$: \ for $v,v'\in \Pi_1$, put
\begin{equation}
\begin{array}{lll}
&&\quad v'\le_{\tilde d}v \\
 &\overset{\textrm{def}}{\Leftrightarrow}&  
\text{the conditions a), b) and c) in Assertion 4.2 iii) hold.}
\end{array}\!\!\!\!\!
\end{equation}

\noindent
In other words,  {\it there is no order relation $\le_{\tilde d}$ between 
$v, v'\!\in\! \Pi_1$ if $\Gamma_{d,v}\!\cap\Gamma_{d,v'}\!=\!\emptyset$, 
otherwise the order relation $\le_{\tilde d}$ agrees with $\le_d$. }
\end{definition}

\begin{assertion}
 Let  $\tilde d\!\in\! \widetilde{Ord}(\Pi_1)$ be given. 
For any 
$v\!\in\!\Pi_1$,
 the set of predecessors 
$\{w\!\in\!\Pi_1\mid w\!<_{\tilde d}\!v\}$ 
is totally ordered by $\le_{\tilde d}$.
\end{assertion}
\begin{proof}
  Suppose $w_i<_{\tilde d}v$  ($i=1,2$).
This means $v\in \Gamma_{d,w_i}$ (Assertion 4.2 iii) c)), and hence 
$\Gamma_{d,w_1}\cap \Gamma_{d,w_2}\not=\emptyset$. Then, Assertion 4.2 i) 
implies that either $w_1\!\le_{\tilde d}\!w_2$ or $w_1\!\ge_{\tilde d}\!w_2$. 
\end{proof}

We obtain the following characterization of 
the partial ordering $<_{\tilde d}$.
\begin{assertion}
 For two orderings $d,d'\in Ord(\Pi_1)$, the following two conditions are 
equivalent.

a)  One has the equality $\Gamma_{d,v}=\Gamma_{d',v}$ for all $v\in \Pi_1$,
i.e.\ $\tilde d=\tilde d'$.

b)  The partial orderings $\le_{\tilde d}$ and $\le_{\tilde d'}$ 
on the set $\Pi_1$ coincide.
\end{assertion}
\begin{proof} 
We have only to show that the partial ordering $<_{\tilde d}$ 
determines the set $\Gamma_{d,v}$ for $v\!\in\!\Pi_1$. 
First, we note that the set $\Gamma_{d,v}$ 
is given by 
\[
\Gamma_{d,v}= \big(\Gamma_{d,v}\cap \Pi_1\big) \cup \bigcup_{w\in \Gamma_{d,v}\cap \Pi_1}
Nbd(w),
\]
where a neighborhood of a point $w\in \Gamma$ is defined by 
\begin{equation}
Nbd(w):=\{u\!\in\! \Pi\mid ^\exists \overline{wu}\in Edge(\Gamma)\}.
\end{equation}
\medskip
\noindent
({\it Proof.} The inclusion $\subset$ follows from the connectivity of 
$\Gamma_{d,v}$. The opposite inclusion $\Gamma_{d,v}\!\supset\! Nbd(w)$ for 
$w\!\in\! \Gamma_{d,v}\!\cap\! \Pi_1$ follows also from the  connectivity of 
$\Gamma_{d,v}$.$\Box$). 
Here, we note that the RHS is determined only from the tree $\Gamma$ and the set 
$\Gamma_{d,v}\cap \Pi_1$, but no ordering is involved.
On the other hand, due to Assertion 4.2  iii) c), one has  
\[
\Gamma_{d,v}\cap \Pi_1:=\{w\in\Pi_1\mid v\le_{\tilde d}w\}.
\]
Thus, $\Gamma_{d,v}$, as a set, is determined by the partial ordering $\le_{\tilde d}$.  
\end{proof}
Due to Assertion 4.4, from now on, we shall identify the equivalent class 
$\tilde d$ with the partial ordering $\le_{\tilde d}$ on $\Pi_1$. 
Thus, {\it the set $\widetilde{Ord}(\Pi_1)$ is regarded as the 
set of certain partial orderings on $\Pi_1$} (see {\it Remark} 6). 
Finally, we give an explicit description of the graph $\Gamma_{d,v}$ 
from   $\le_{\tilde d}$.

\begin{assertion}
Let a linear ordering for $d\in Ord(\Pi_1)$ be given. One has  
the equality $\Gamma_{d,v}=\Gamma_{\tilde d,v}$ for any $v\in \Pi_1$, where 
\begin{equation}
\begin{array}{lll}
\Gamma_{\tilde d,v}&:=&\text{the connected component of 
$\Gamma\setminus \{w\!\in\! \Pi_1\mid w\!<_{\tilde d}\!v\}$}\\
&& \text{ containing $v$. }
\end{array}
\end{equation}
\end{assertion}
\begin{proof}
Since the total ordering $d$ is a refinement of $\tilde d$,
one has the inclusion  $\Gamma_{d,v}\subset\Gamma_{\tilde d,v}$.
To show the opposite inclusion, it is sufficient to show that if 
$v\not\le_{\tilde d}w$, then $w$ does not belong to $\Gamma_{\tilde d,v}$.

We may assume $w\not\le_{\tilde d}v$, otherwise $w\not\in \Gamma_{\tilde d,v}$ 
is trivial.
Consider the totally ordered set
$\{u\in\Pi_1\mid u<_{\tilde d}w,\ u<_{\tilde d}v\}$ (c.f.\ Fact iv)). 
Since it contains at least the minimum element of $d$ and non-empty, 
it contains the unique maximal element, say $u_0$.
By definition, $v$ and $w$ belong to the connected component 
$\Gamma_{d,u_0}\subset\Gamma_{\tilde d,u_0}$. However, since they 
belong to different components of $\Gamma_{d,u_0}\setminus \{u_0\}$,
they also belong to different components of 
$\Gamma_{\tilde d,u_0}\setminus \{u_0\}$
(since $\Gamma_{\tilde d,u_0}$ is a tree).
\end{proof}

Let us return to the proof of the Theorem.
The formula (19) is shown by induction on $l=\#\Gamma$. 
We first prepare an induction formula.

Let $\Gamma$ be a tree. For a given principal decomposition 
$\{\Pi_1,\Pi_2\}$ and an  
attached principal orientation $o_{\Pi_1,\Pi_2}$, 
we enumerate the set $\Sigma(o_{\Pi_1,\Pi_2})$.

By definition, for any total ordering $d\in \Sigma(o_{\Pi_1,\Pi_2})$,
the smallest element belongs to $\Pi_1$. Therefore, we have a decomposition:
\[
\Sigma(o_{\Pi_1,\Pi_2})=\coprod_{v\in\Pi_1} \Sigma(o_{\Pi_1,\Pi_2},v<)
\]
where 
$\Sigma(o_{\Pi_1,\Pi_2},v<):=\{d\in \Sigma(o_{\Pi_1,\Pi_2}) \mid 
v \text{ is the smallest element in $d$.}\}$. 
Put 
$\sigma(o_{\Pi_1,\Pi_2},v<):=\#\Sigma(o_{\Pi_1,\Pi_2},v<)$ so that one has 
\[
\sigma(\Gamma)=\sigma(o_{\Pi_1,\Pi_2})
=\sum_{v\in \Pi_1}\sigma(o_{\Pi_1,\Pi_2},v<).
\]
For $w\in Nbd(v)$,
let us denote by $\Gamma_{vw}$ the connected component of 
$\Gamma\setminus\{v\}$
containing $w$. One has the decomposition 
$\Gamma\setminus\{v\}=\coprod_{w\in Nbd(v)} \Gamma_{vw}$. 

Applying (10) in Assertion 2.1 for $\al\!=\!v$ and 
$r\!=\!l\!-\!1\!=\!\sum_{w\in Nbd(v)}\!r_w$, 
$r_w:=\#\Gamma_{vw}$, we obtain:
$ \sigma(o_{\Pi_1,\Pi2},v<)=(l-1)! \prod_{w\in Nbd(v)}\frac{\sigma(\Gamma_{vw})}{(\#\Gamma_{vw})!}. $

\noindent
Summing over all vertices $v\!\in\! \Pi_1$, 
we obtain the induction formula:
\begin{equation}
\frac{\sigma(\Gamma)}{(\#\Gamma)!}=\frac{1}{l} \sum_{v\in\Pi_1}\prod_{w\in Nbd(v)}\frac{\sigma(\Gamma_{vw})}{(\#\Gamma_{vw})!}.
\end{equation}

By the induction hypothesis, for any $v\in \Pi_1$ and $w\in Nbd(v)$, 
we have already the formula for $\Gamma_{vw}$:
\[
\frac{\sigma(\Gamma_{vw})}{(\#\Gamma_{vw})!}=
 \sum_{\tilde d_w
\in \widetilde{Ord}(\Gamma_{vw}\cap\Pi_1)}\frac{1}{\prod_{u\in \Gamma_{vw}\cap\Pi_1} 
\#(\Gamma_{vw})_{\tilde d_w,u}}.
\leqno{*)}\] 
The substitution of  $*)$ into RHS of (23) gives a formula
summing the terms: $\frac{1}{l}\prod_{w\in Nbd(v)}\frac{1}{\prod_{u\in \Gamma_{vw}\cap\Pi_1} 
\#(\Gamma_{vw})_{\tilde d_w,u}}$, where the summation index 
$v\times \{\tilde d_w\}_{w\in Nbd(v)}$ runs over the set 
$\bigcup_{v\in\Pi_1}\big(v\times \prod_{w\in Nbd(v)}( \widetilde{Ord}(\Gamma_{vw}\cap\Pi_1))\big).$

To the index $v\times \{\tilde d_w\}_{w\in Nbd(v)}$, we attach 
the partial ordering $\tilde d$ of the set $\Pi_1$
defined by the rule a) $v$ is the smallest element,
b) $\tilde d$ agrees with $\tilde d_w$ on the set $\Gamma_{vw}\cap \Pi_1$
for $w\in Nbd(v)$, and c) there is no order relation between 
$\Gamma_{vw}\cap \Pi_1$ and $\Gamma_{vw'}\cap\Pi_1$  
for different $w,w'\!\in\! Nbd(v)$.

The correspondence $v\times \{\tilde d_w\}_{w\in Nbd(v)}\mapsto \tilde d$ 
gives a bijection:
\[
\bigcup_{v\in\Pi_1}\big(v\times \prod_{w\in Nbd(v)}( \widetilde{Ord}(\Gamma_{vw}\cap\Pi_1))\big) 
\ \simeq \ \widetilde{Ord}(\Pi_1),
\]
where the opposite correspondence is given by the restriction map.

On the other hand, the term 
$\frac{1}{l}\prod_{w\in Nbd(v)}\frac{1}{\prod_{u\in \Gamma_{vw}\cap\Pi_1} 
\#(\Gamma_{vw})_{\tilde d_w,u}}$ 
for the index $v\times \{\tilde d_w\}_{w\in Nbd(v)}$
coincides with the term 
$\frac{1}{\prod_{u\in \Pi_1} \#\Gamma_{\tilde d,u}}$ in (19) given by 
the corresponding partial ordering $\tilde d$.
This means that the substitution of $*)$ into RHS of (23)
gives RHS of the formula (19). 

This completes the proof of the Theorem 4.1.
\end{proof}

 \begin{remark}
In the next \S5, the set $\widetilde{Ord}(\Pi_1)$ plays again quite an important 
role, where we regard $\widetilde{Ord}(\Pi_1)$ as the set of reduced oriented 
graphs 
(see Theorem 5.1 and its corollaries).
\end{remark}

\section{ Block decomposition of the principal $\Gamma$-cone}

As the third main result of the present paper, we decompose  
a principal $\Gamma$-cone $E_{o_{\Pi_1,\Pi_2}}$ into blocks,
where each block is a simplicial cone associated to a {\it rooted tree} 
and is characterized combinatorially. 
The number of chambers in a block is  
given by the hook length formula. Thus, 
we obtain an alternative but intrinsic proof of the formula (19).

\begin{definition} 
A reduced oriented graph $(\Gamma,o)$ is called a {\it rooted tree} if 

i)  There exists a unique minimal vertex $v_{o}\in\Gamma$ 
with respect to $o$.

ii)  Any vertex ($\not=v_{o}$) of $\Gamma$  has a 
unique immediate predecessor.

\noindent
The smallest vertex $v_o$ is called the {\it root} of $(\Gamma,o)$.
{\it The definition implies
that $\Gamma$ is a tree}. Conversely, a pair of 
a tree $\Gamma$ and a vertex $v_o$ of $\Gamma$ determines a unique  
rooted tree structure having $v_o$ as its root.
\end{definition}

We return to the setting of \S4: $\Gamma$ is a tree on 
$\Pi$, $\Pi=\Pi_1\amalg \Pi_2$ is a principal decomposition (13), 
and $\widetilde{Ord}(\Pi_1)$ is the set of certain partial 
orderings on $\Pi_1$ (recall \S4 Definition (20),
the equivalence of a) and b) in Assertion 4.4 and the paragraph following to 
Assertion 4.4).  
We, further, identify the partial ordering 
$\tilde d\!\in\! \widetilde{Ord}(\Pi_1)$ with its reduced 
oriented graph $\tilde d\!=\! \tilde d_{red}$ 
on the set $\Pi_1$ (= the Hasse diagram,
see the paragraph before {\it Remark} 1. in \S1).
That is: $\widetilde{Ord}(\Pi_1)$ is regarded as a set of certain 
reduced oriented graph structures on $\Pi_1$.
The followings are reformulations of what we have shown in \S4.

\medskip
\noindent
{\bf Fact.} 
{\bf a)} {\it Any element $\tilde d\!\in\!\widetilde{Ord}(\Pi_1)$ 
is a rooted tree structure on $\Pi_1$. 
} 

\noindent
{\bf b)} {\it For any total ordering $d\!\in\! Ord(\Pi_1)$, there exists a unique
$\tilde d\!\in\! \widetilde{Ord}(\Pi_1)$ such that $d$ is 
a linear extension of $\tilde d$.}

\noindent
{\bf c)} 
{\it The system  
$\{E_{\tilde d}\}_{\tilde d \in \widetilde{Ord}(\Pi_1)}$ 
is a simplicial cone decomposition of $V_{\Pi_1}$.}

\begin{proof} a) 
Let $\tilde d$ be the equivalence class of $d\!\in\! Ord(\Pi_1)$. 
The smallest element, say $v_d$, of $d$ is also the smallest 
with respect to $\tilde d$, since $\Gamma_{v_d}\!=\!\Gamma$ contains all $\Pi_1$ 
(c.f.\ Assertion 4.2, ii) c)). 
Then, the uniqueness of the predecessor for an element $v\not =v_d$ follows from Assertion 4.3. 

b) This follows from the definition in \S4 of $\widetilde{Ord}(\Pi_1)$,
where any total ordering $d$ belongs to the unique equivalence class $\tilde d$.

c) This follows from a), b) and Assertion 3.1.
\end{proof}


In the following Theorem 5.1 and its corollary, 
we lift (in a suitable sense) the above Facts to the principal $\Gamma$-cone $E_{\Pi_1,\Pi_2}$.

First, let us define the lifting $\tilde{o d}$ of 
$\tilde d\in \widetilde{Ord}(\Pi_1)$ 
by 
\begin{equation}
\tilde d\mapsto \tilde{od}:=(o_{\Pi_1,\Pi_2}\cup \tilde d)_{red},
\end{equation}
where $(o_{\Pi_1,\Pi_2}\cup \tilde d)_{red}$ is the reduction, 
i.e.\ Hasse diagram (recall \S1), 
of the oriented graphs $o_{\Pi_1,\Pi_2}\!\cup\! \tilde d$ 
obtained by the union of oriented edges of $o_{\Pi_1,\Pi_2}$ and $\tilde d$
(here, the union is acyclic, since  i) any element of 
$\Pi_2$ cannot be a part of an oriented cycle since it 
is maximal with respect to $o_{\Pi_1,\Pi_2}\cup \tilde d$, and  
ii) the part $\tilde d$ on $\Pi_1$ is a tree without a cycle, Fat a)).

The following theorem gives a characterization of the elements 
$\tilde d\in \widetilde{Ord}(\Pi_1)$ and their liftings $\tilde{od}$. 
Then the corollary shows that the cones $E_{\tilde{od}}$  for the liftings give arise a 
simplicial decomposition of the principal cone $E_{\Pi_1,\Pi_2}$,
which we shall call the {\it block decomposition} of $E_{\Pi_1\Pi_2}$.

\begin{theorem} 
Let $\Gamma$ be a tree on $\Pi$, and let $\Pi=\Pi_1\amalg\Pi_2$ {\rm (13)}
be one of its principal decompositions.

{\bf 1.} The set $\widetilde{Ord}(\Pi_1)$ of rooted trees on $\Pi_1$
is characterized as follows:
\[
\ \ \begin{array}{lll}
\widetilde{Ord}(\Pi_1)=&\!\!\! \{\ \tilde d\ \mid 
\text{ {\rm i)}  $\tilde d$ is a rooted tree structure on $\Pi_1$,} \\
&\text{ {\rm ii)} $\tilde d|_{Nbd(\beta)}$ is totally ordered for any $\beta\in\Pi_2$,}\\
& \text{{\rm iii)} $\tilde d$ is minimal with respect to {\rm i)} and {\rm ii)},}\\
& \text{\ \ i.e.\ if $\tilde f$ satisfies {\rm i), ii)} and $\tilde f\!\le\!\tilde d$ 
then $\tilde f\!=\!\tilde d$.}
\}
\end{array}
\vspace{-0.1cm}
\]
Here, we recall {\rm (21)} for a definition of a neighborhood $Nbd(\beta)$ of $\beta$.

{\bf 2.} Let us introduce a set of rooted tree structures on $\Pi$:  
\[
\begin{array}{lll}
\widetilde{Ord}_{\Pi_1,\Pi_2}(\Pi):=&\!\!\! \{\ \tilde e\ \mid 
\text{ {\rm i)}  $\tilde e$ is a rooted tree structure on $\Pi$,} \\

&\text{ {\rm ii)} $o_{\Pi_1,\Pi_2}\le \tilde e$, i.e.\ if $x\!<_{o_{\Pi_1,\Pi_2}}\!\!y$ 
then $x\!<_{\tilde e}\!y$},\\
& \text{{\rm iii)} the set $\Pi_2$ is totally disordered with repect to $\tilde e$,}\\
& \text{{\rm iv)} $\tilde e$ is minimal with respect to {\rm i), ii)} 
and {\rm iii)},}\\
& \text{ \ i.e.\ if $\tilde f$ satisfies {\rm i), ii), iii)} and $\tilde f\!\le\!\tilde e$ 
then $\tilde f\!=\!\tilde e$.}
\}
\end{array}
\]
Then, the lifting 
$\tilde d \mapsto \tilde{od}$ {\rm (24)} induces the bijection
\begin{equation}
\widetilde{Ord}(\Pi_1)\ \simeq\ \widetilde{Ord}_{\Pi_1,\Pi_2}(\Pi).
\end{equation}

{\bf 3.} The set  $\Sigma(o_{\Pi_1,\Pi_2})$ of chambers in $E_{\Pi_1,\Pi_2}$  decomposes into a union: 
\begin{equation}
\Sigma(o_{\Pi_1,\Pi_2})\ =\  \amalg_{\tilde{e}\in \widetilde{Ord}_{\Pi_1,\Pi_2}(\Pi)}\Sigma(\tilde{e}).
\end{equation}
\end{theorem}
\begin{definition}
We call the simplicial cone $E_{\tilde{e}}$ associated to a rooted tree 
$\tilde e\in \widetilde{Ord}_{\Pi_1,\Pi_2}(\Pi)$
a {\it block}.
\end{definition}
The decomposition (26) can be paraphrased in terms 
of the block decomposition of the principal $\Gamma$-cone as follows.
\begin{corollary} (Block decomposition of the principal cone)

The principal cone $E_{\Pi_1,\Pi_2}$ decomposes into 
a union of the blocks: 
\begin{equation}
\overline E_{\Pi_1,\Pi_2} \ =\ \coprod_{\tilde e\in\widetilde{Ord}_{\Pi_1,\Pi_2}(\Pi)} 
\overline E_{\tilde{e}} \ .
\end{equation}


\end{corollary}
\begin{proof} 
We prove only {\bf 2}, and then {\bf 1.}\ is its byproduct.
The proof is slightly involved and is    
divided into steps a.i)-ii),b.i)-ii),c),d) and e).

a.i) {\it The lifting $\tilde{od}$ 
for $\tilde d\!\in\! \widetilde{Ord}(\Pi_1)$ satisfies i), ii) and iii) of {\bf 2}.}

a.ii) $*)$ {\it For any $\beta\!\in\! \Pi_2$ and for any 
 $\tilde d\!\in\! \widetilde{Ord}(\Pi_1)$, 
the set $Nbd(\beta)\subset\Pi_1$
is totally ordered with respect to the partial ordering $\tilde d$.}

$**)$ {\it The restriction $\tilde {od}|_{\Pi_1}$ of the 
oriented graph $\tilde {od}$ to $\Pi_1$ is equal to $\tilde d$. 
In particular, this implies that the correspondence {\rm (24)} is injective.}

\vspace{-0.1cm}

\begin{proof} 
a.i) Those ii) and iii) for $\tilde{od}$ follow immediately from the definition. 
The  i) for $\tilde{od}$ is a consequence of a.ii) and b.ii).

a.ii)$*)$ For $\beta\!\in\!\Pi_2$, 
consider  $\al_1,\al_2\!\in\! Nbd(\beta)\!\subset\! \Pi_1$ with $\al_1\!\not=\!\al_2$.
Since $\Gamma_{\tilde d,\al_i}$ (recall (22)) contains $Nbd(\al_i)$,
one has $\beta\!\in\! \Gamma_{\tilde d,\al_i}$ for $i\!=\!1,2$.
That is $\Gamma_{\tilde d,\al_1}\!\cap\! \Gamma_{\tilde d,\al_2}\!\not=\!\emptyset$.
Then, due to Assertion 4.2, i) and ii), either 
$\Gamma_{\tilde d,\al_1}\!\supset\! \Gamma_{\tilde d,\al_2}$ or 
$\Gamma_{\tilde d,\al_1}\!\subset\! \Gamma_{\tilde d,\al_2}$ occurs, 
and, hence by the definition (20), one has either 
$\al_1\!<_{\tilde d}\!\al_2$ or $\al_1\!>_{\tilde d}\!\al_2$. 
Thus, $*)$ is shown.

a.ii)$**)$ By definition, $\tilde {od}|_{\Pi_1}\!\le\!\tilde d$. 
If an edge $\overline{\al\al'}$ with $\al\!<_{\tilde d}\!\al'$ in $\tilde d$ 
is removable in $\tilde {od}$, there exists a path 
$\al_0\!\!=\!\!\al\!\!<_{\tilde {od}}\!\al_1\!\!<_{\tilde {od}}\!\cdots\!\!<_{\tilde {od}}\!\al_n\!\!=\!\!\al'$ 
in $\tilde {od}$ for $n\!\ge\!2$. Since $\tilde d$ is reduced, there exits 
$1\!\!\le\!\! i\!\!<\!\!n$ such that $\al_i\!\!\in\!\!\Pi_2$. 
By taking largest such $i$, 
we find an edge $\al_i\!\!<_{\tilde {od}}\!\al_{i+1}$ with $\al_i\!\!\in\!\!\Pi_2$ and 
$\al_{i+1}\!\!\in\!\!\Pi_1$.\ This is impossible 
since $\tilde {od}\!\!\le\!\!(o_{\Pi_1\Pi_2}\!\cup\!\tilde d)$.\ 
Thus, $**)$ is shown.
\end{proof}
\vspace{-0.1cm}

b.i) {\it Let $\tilde e$ be an oriented graph on $\Pi$
satisfying ii) and iii) of {\bf 2}. Then,
there does not exists a pair $\al\!\in\!\Pi_1$ and $\beta\!\in\!\Pi_2$ such that 
$\beta\!<_{\tilde e}\!\al$. }

b.ii) 
{\it Under the same assumption, $\tilde e$ is a rooted tree   
(i.e.\ $\tilde e$ satisfies {\rm i)}), 
if and only if one has 
$*)$ $\tilde e|Nbd(\beta)$ is totally ordered for all $\beta\in \Pi_2$ 
and $**)$ the restriction $\tilde e|_{\Pi_1}$ of the 
oriented graph $\tilde e$ is a rooted tree on $\Pi_1$.
In this case, the restriction to $\Pi_1$ of the ordering $\le_{\tilde e}$ on $\Pi$
coincides with the ordering on $\Pi_1$ generated by the oriented graph 
$\tilde e|_{\Pi_1}$ on $\Pi_1$.
} 
\vspace{-0.1cm}

\begin{proof}
b.i) For any $\al\in\Pi_1$, there exists at least 
one $\gamma\!\in\! \Pi_2$ such that $\al\!\!<_{o_{\Pi_1,\Pi_2}}\!\!\gamma$.
If the contrary $\beta\!<_{\tilde e}\!\al$ holds,  using ii) and the transitive 
closure,
one obtains $\beta\!\!<_{\tilde e}\!\!\gamma$, which contradicts to iii).

b.ii)$*)$ If there were no order relation among $\al_1,\al_2\!\in\! Nbd(\beta)$, there 
exists some $\gamma\!\in\! \Pi$ such that 
$\al_1,\al_2\!<_{\tilde e}\!\gamma\!\le_{\tilde e}\!\beta$ and $\gamma$ has at least two 
immediate predecessor. A contradiction!
$**)$ By b.i), the set of predecessors for any $\al\!\in\! \Pi_1$ in $\Pi$ 
is equal to that in $\Pi_1$. In particular, the imediate predecessor $\al'$
is in $\Pi_1$ so that the edge $\overline{\al'\al}$ belongs to 
$\tilde e|_{\Pi_1}$.  

Conversely, if $Nbd(\beta)$ is totally ordered, then $\beta\!\in\!\Pi_2$
is connected with only the maximal element of $Nbd(\beta)\!\subset\! \Pi_1$. 
This, together with $**)$, implies that $\tilde e $ is a rooted tree.
\vspace{-0.2cm}
\end{proof}

c) {\it The lifting $\tilde{od}$ for $\tilde d\!\in\! \widetilde{Ord}(\Pi_1)$ satisfies iv) of Theorem {\rm 5.2}.}
\vspace{-0.1cm}

\begin{proof} 
Assume the contrary: $\exists\tilde f\!\le\! \tilde{od}$ and $\tilde f\!\not=\! \tilde{od}$. 

We, first, show that there exists a pair $\al,\beta\!\in\! \Pi_1$ 
such that $\al\!<_{\tilde{od}}\!\beta$ 
and $\al\!\not<_{\tilde f}\!\beta$,
$\al\!\not>_{\tilde f}\!\beta$. 
({\it Proof.} By the assumption, there exist at least a pair $\al,\beta\!\in\!\Pi$ 
with  $\al\!<_{\tilde{od}}\!\beta$ and $\al\!\not<_{\tilde f}\!\beta$.
Due to iii) and b.i) for $\tilde{od}$, one has $\al\!\not\in\! 
\Pi_2$. 
 If $\al\!\in\! \Pi_1$ and $\beta\!\in\!\Pi_2$ then by
construction of $\tilde{od}$ 
there exists $\gamma\!\in\! \Pi_1$ such that $\al\!<_{\tilde d}\!\gamma$ and
$\gamma\!<_{o_{\Pi_1,\Pi_2}}\!\beta$. Then, one necessary has 
$\al\!\not<_{\tilde f}\!\gamma$,
since otherwise one has $\al\!<_{\tilde f}\!\gamma\!<_{o_{\Pi_1,\Pi_2}}\!\beta$,
 and by ii) on $\tilde f$, 
one has the contradiction $\al\!<_{\tilde f}\!\beta$. Then,
replace $\beta$ by $\gamma$.) 

We may assume that $\al$ is the immediate 
predecessor of $\beta$ with respect to $\tilde d$, since one of the edges of 
the path from $\al$ to $\beta$ satisfies the condition. 
We may assume further that the pair $\al,\beta\!\in\! \Pi_1$ is maximal in the sense 
that
for any $\beta'\!\in\! \Gamma_{\tilde d,\beta}\cap \Pi_1$ 
(i.e.\ $\beta\!\le_{\tilde d}\!\beta'$) 
one has $\beta\!\le_{\tilde f}\!\beta'$. 

Since $\al$ is the immediate predecessor of $\beta$ with respect to $\tilde d$,
$\al$ is connected with a single point, say $\gamma\in \Pi_2$, of 
$\Gamma_{\tilde d,\beta}$. Since $\Gamma_{\tilde d,\beta}$ is a connected tree 
containing $\beta$, there exists a path in $\Gamma_{\tilde d,\beta}$ connecting 
$\gamma$ and $\beta$. Let $\beta'\in \Pi_1$ be the element next to $\gamma$. 
Since $\beta'\in \Gamma_{\tilde d,\beta}\cap\Pi_1$, 
by the assumption of the maximality of the pair $\al,\beta$, we have 
$\beta\le_{\tilde f}\beta'$. 
On the other hand, due to ii) for $\tilde f$, we have also 
$\al\le_{\tilde f}\gamma$ and $\beta'\le_{\tilde f}\gamma$. Thus we obtain two 
monotonically increasing pathes attached to the two sequences
$v_{\tilde{f}}\le_{\tilde{f}}\al\le_{\tilde{f}}\gamma$ and
$v_{\tilde{f}}\le_{\tilde{f}}\beta\le_{\tilde{f}}\beta'\le_{\tilde{f}}\gamma$. 
They are different, since there is no order relations between $\al$ and $\beta$ 
with respect to $\tilde f$.
A contradiction to that $\tilde f$ is a tree!
 Thus, iv) for $\tilde{od}$ is shown.
\end{proof}

d) {\it Any element $\tilde e\in \widetilde{Ord}_{\Pi_1,\Pi_2}(\Pi)$ 
is a lifting $\tilde{od}$ for some $\tilde d\!\in\! \widetilde{Ord}(\Pi_1)$.}
\begin{proof}

 Consider any linear extension $d\!\in\! Ord(\Pi_1)$ of $\tilde e|_{\Pi_1}$ 
(i.e.\ $e|_{\Pi_1}\!\le\! d$) and its  equivalence class 
$\tilde d\!\in\! \widetilde{Ord}(\Pi_1)$ (\S4 Definition). 
If we have shown 
$\tilde d\!\le\!\tilde e|_{\Pi_1}$ (shown in the following e)), 
then we have done the proof since
we have 
$\tilde{od}\!=\!(o_{\Pi_1,\Pi_2}\!\cup\!\tilde d)_{red}\!\le\! 
(o_{\Pi_1,\Pi_2}\!\cup\!\tilde e|_{\Pi_1})_{red}\!\le\! 
\tilde e$ and the minimality iv) for $\tilde{e}$ implies 
$\tilde{od}\!=\!\tilde e$. 
Thus, the surjectivity of (25) is shown.
\end{proof}

e) {\it For any $d\in Ord(\Pi_1)$ with $\tilde e|_{\Pi_1}\le d$, one has 
$\tilde d\!\le\!\tilde e|_{\Pi_1}$.}

\begin{proof} For 
any edge $\overline{\al\beta}$ in $\tilde d$ on $\Pi_1$ with 
$\al\!<_{\tilde d}\!\beta$, 
we have to show the order relation $\al\!<_{\tilde e}\!\beta$.
We show this by the ``inverse induction'' in the sense that we prove  
it by assuming the innequality $\al'<_{\tilde e}\beta'$ for 
all edges $\overline{\al'\beta'}$  in $\tilde d$ with the inequality: 
$\beta\le_{\tilde d}\al'<_{\tilde d}\beta'$.

Suppose the contrally $\al\not<_{\tilde e}\beta$. Since 
$\al>_{\tilde e}\beta$ is not possible (else, we have $\al>_{d}\beta$, contradicting 
to $\al<_{\tilde d}\beta$), there should be no order relation between 
$\al$ and $\beta$ with respect to $\tilde e$.
The fact $\al<_{\tilde d}\beta$ implies that $\beta\in \Gamma_{d,\al}$ (recall (20)). 
The fact that $\al$ is the immediate predecessor of $\beta$ means, precisely, 
that there exist elements $\beta'\in 
\Gamma_{d,\beta}\cap \Pi_1$ and $\gamma\in \Pi_2$ such that 
$\al,\beta'<_{o_{\Pi_1,\Pi_2}}\gamma$. By assumption ii) on $\tilde e$, 
we have $\al,\beta'<_{\tilde e}\gamma$.
On the other hand, by definition $\beta'\in \Gamma_{d,\beta}$ and hence $\beta\le_{\tilde d}\beta'$. 
Then, by the inverse induction hypothesis, one obtains $\beta\le_{\tilde e}\beta'$.
Thus, we obtain two monotonically increasing pathes in $\tilde e$ connecting 
the bottom element $v_{\tilde e}$ and $\gamma$ 
as follows: $v_{\tilde e}\le_{\tilde e}\al\le_{\tilde e}\gamma$ and 
$v_{\tilde e}\le_{\tilde e}\beta\le_{\tilde e}\beta'\le_{\tilde e}\gamma$.
Since, by the contraly assumption, 
there is no order relation between $\al$ and $\beta$ with respect to $\tilde e$,
whereas each path passes through either $\al$ or $\beta$, the two pathes are different.
This is a contradiction to the fact that $\tilde e$ is a tree. 
This proves $\al<_{\tilde e}\beta$ and, hence, e) is proven.
\end{proof}

3.
For $\tilde e\! \in\! \widetilde{Ord}_{\Pi_1,\Pi_2}(\Pi)$, ii) implies 
$\Sigma(o_{\Pi_1,\Pi_2})\!\supset\!\Sigma(\tilde{e})$. 
The decomposition (26) follows, since for any $d\!\in\! \Sigma(o_{\Pi_1,\Pi_2})$,
there exists a unique $\tilde d\!\in\! \widetilde{Ord}(\Pi_1)$ such that 
$(o_{\Pi_1,\Pi_2}\!\cup \tilde d)\le d$ ({\bf Fact}\! b) and c)). 

Corollary is a rewriting of (26) (c.f.\ 2. of Assertion 1.1).

These complete the proof of Theorem 5.1 and its Corollaries.
\end{proof} 
\begin{remark}
The block decomposition (27) of the principal $\Gamma$-cone 
$E_\Gamma$ depends on a choice of the principal 
orientation $o_{\Pi_1,\Pi_2}$ 
(see Example).
\end{remark}

As an application of the block decomposition (26) ((27)) of the principal cone, 
let us give 
an alternative proof of the formula (19).
This is achieved by two steps.
The first step is to recall the well known 
hook length formula of Knuth enumerating the chambers 
in a $\Gamma$-cone for a rooted tree (it is also an immediate 
consequence of the decomposition formula (10) with $r\!=\!\#\Gamma\!-\!1$ and $\al\!=$
the root of the tree).

\begin{lemma} {\rm (Knuth [K2,p70])}
Let $(\Gamma,o)$ be a rooted tree.
Then 
\begin{equation}
\sigma(o)=\frac{(\#\Gamma)!}{\prod_{v\in \Pi} 
\#\Gamma_{o ,v}}
\vspace{-0.5cm}
\end{equation}
where
\begin{equation}
\begin{array}{lll}
\Gamma_{o,v}\!&\!:=\!&\!\text{the connected component of $\Gamma_{o}
\setminus \{w\in\Pi\mid w<_{o}v\}$}\\ 
&&\text{ containing $v$.}
\end{array}\!\!\!
\end{equation}
\end{lemma}
\begin{note}
The underlying 
oriented graph structure in $\Gamma_{\tilde d,v}$ (22) is the principal 
orientation $o_{\Pi_1,\Pi_2}$ and that for  $\Gamma_{o,v}$ (29)
is the rooted tree $o$.
They are, in a sense, the most contrasting orientations.
However, we show in the following a ``numerical coincidence'' of them.
\end{note}

The second step of the alternative proof of (19) is as follows.
Apply (28) to $\sigma(\tilde{od})$ 
to count the number of chambers in $E_{\tilde{od}}$. 
Comparing (19) and (25),
let us show the equality:
\begin{equation}
\quad\quad\quad\quad\quad
\frac{(\#\Gamma)!}{\prod_{v\in \Pi}\#\Gamma_{\tilde{od},v}}
= \frac{(\#\Gamma)!}{\prod_{v\in \Pi_1} 
\#\Gamma_{\tilde{d},v}} 
\end{equation}
for $ \tilde d \in \widetilde{Ord}(\Pi_1).$
We note that the region of the running index $v$ in 
%
LHS of (30) can 
be shrunken from $\Pi$ to $\Pi_1\!=\!\Pi\!\setminus\!\Pi_2$, 
since for $v\in \Pi_2$ one has $\#\Gamma_{\tilde{o d},v}\!=\!1$
because of the fact that $v$ is maximal with respect to 
$o_{\Pi_1,\Pi_2}\!\cup\! \tilde d$.
Therefore, we have only to show the formula
\begin{equation}
\#\Gamma_{\tilde{d},v}=\#\Gamma_{\tilde{o d},v}
\end{equation}
for $v\in \Pi_1$. We show that the vertex sets 
$|\Gamma_{\tilde{d},v}|$ and $|\Gamma_{\tilde{o d},v}|$ coincide 
(even though the graph structures are quite different, see Example).

Note the inclusion relation:\! 
$o_{\Pi_1,\Pi_2}\! \subset\! (o_{\Pi_1,\Pi_2}\!\cup\tilde d)\! \supset\! \tilde{o d}$ 
among oriented graphs and the equality among the vertex sets:\! 
 $A\!:=\!\{w\!\in\!\Pi_1\!\mid w\!<_{\tilde d}v\}\!=\!
\{w\!\in\!\Pi\!\mid w\!<_{o_{\Pi_1,\Pi_2}\! \cup\!\tilde d}v\}\!=\!\{w\!\in\!\Pi\!\mid w\!<_{\tilde{o d}}v\}$.
Thus, one has the relation: $\Gamma_{\tilde d,v}\!\subset\! \Gamma_{o_{\Pi_1,\Pi_2}\cup\tilde d,v}\!\supset\! \Gamma_{\tilde{o d},v}$
among the connected components of the complements of $A$ containing $v$.
The sets $|\Gamma_{\tilde d,v}|$ and $|\Gamma_{o_{\Pi_1,\Pi_2}\cup\tilde d,v}|$ coincide, since,
if $v\!<_{\tilde d}\!w$ for $w\!\in\!\Pi_1$ then $w\!\in\! \Gamma_{\tilde{d},v}$. 
The sets $|\Gamma_{o_{\Pi_1,\Pi_2}\cup\tilde d,v}|$ and 
$|\Gamma_{\tilde{o d},v}|$ coincide, since, 
if an element $w\!\in\! \Pi_2$ is connected with $v$ in 
$\Gamma_{o_{\Pi_1,\Pi_2} \cup\tilde{d},v}$ then $w$ is connected with $v$ by $\tilde{od}$.

This completes the proof of the equality (31), 
and hence the alternative proof for (19) is completed.

\medskip
\noindent
{\bf Example.} \
We illustrate the two block decompositions of type $A_7$ and the 
formula (19) for the two principal orientations on $\Gamma(A_7)$. 

\medskip
\noindent 
{\bf I.} \ 
Let the principal decomposition (13) of $\Gamma(A_7)$ 
and the associated principal orientation $o_{A_7}:=o_{\Pi_1,\Pi_2}$ be given by 
\[
\begin{array}{llllllllll}
\vspace{-0.2cm}
&\Pi_2:& \circ\!\! &\!\!&\!\!\!\!\circ\!\!&&\!\!\!\!\circ&&\!\!\!\!\!\circ& \\
\vspace{-0.2cm}
o_{A_7}:&& \ \ \nwarrow&\nearrow&\!\nwarrow&\nearrow&\!\nwarrow&\nearrow   \\
&\Pi_1:&&\!\!\!\!\circ &&\!\!\!\!\circ&&\!\!\!\!\circ&& .\\
\end{array}
\vspace{-0.1cm}
\]

There are 5 partial orderings $\tilde d\in  \widetilde{Ord}(\Pi_1)$ and,
accordingly,  
5 associated rooted trees $\tilde{od}$. Three of them are illustrated as follows. Here, the black colored vertex is the root of the tree:
\vspace{-0.1cm}
\[
\begin{array}{llllllllllllllllllll}
\vspace{-0.3cm}
& \circ &&\!\!\!\!\circ&\!\!&\!\!\!\!\!\circ\!&\!\!&\!\!\!\!\!\circ&\qquad &\circ\!\! &&\!\!\!\!\circ\!&&\!\!\!\!\circ\!&&\!\!\!\!\!\!\circ&\\
\vspace{-0.2cm}
\!\!\!o_{A_7}\!\cup\tilde d_1:\!\!\!& \ \ \nwarrow&\nearrow&\!\nwarrow&\nearrow&\!\!\nwarrow&\nearrow &&\!\!\!\overset{red}{\longrightarrow}\  \tilde{o d_1}:\!\!\!  & \ \ \nwarrow&&\!\nwarrow&&\!\nwarrow&\nearrow & \\
&&\!\!\!\!\bullet &\!\!\!\!\!\!\!\!\longrightarrow\!\!\!\!&\!\!\!\!\circ\!&\!\!\!\!\!\!\!\longrightarrow\!\!&\!\!\!\!\circ&& &&\!\!\!\!\bullet &\!\!\!\!\!\longrightarrow&\!\!\!\!\circ&\!\!\!\!\!\!\!\longrightarrow&\!\!\!\!\circ&& \\
\end{array}
\]
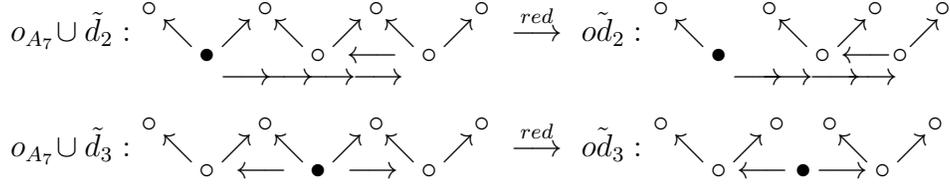
\vspace{-0.1cm}
\[
\begin{array}{llllllllllllllllllll}
\vspace{-0.3cm}
& \circ &&\!\!\!\!\circ&&\!\!\!\!\circ&&\!\!\!\!\!\circ&\qquad &\circ\!\!&&\!\!\!\circ&&\!\!\!\!\!\circ&&\!\!\!\!\!\!\circ&\\
\vspace{-0.2cm}
\!\!\!o_{A_7}\!\cup\tilde d_2:\!\!\!& \ \ \nwarrow&\nearrow&\!\nwarrow&\nearrow&\!\nwarrow&\nearrow &&\!\!\! \overset{red}{\longrightarrow}\  \tilde{o d_2}:\!\!\! & \ \ \nwarrow&&\nwarrow&\nearrow&&\!\!\!\nearrow & \\
\vspace{-0.2cm}
&&\!\!\!\!\bullet &\!\!\!\!\!\!\!\!&\!\!\!\!\circ&\!\!\!\!\!\!\!\!\longleftarrow\!\!\!&\!\!\!\!\circ&& &&\!\!\!\!\bullet &&\!\!\!\!\circ&\!\!\!\!\!\!\!\!\!\!\longleftarrow\!\!\!&\!\!\!\!\!\!\circ&& \\
&&  \longrightarrow\!\!\!\!&\!\!\! \longrightarrow\!\!\!\!&\!\!\!\!\!  \longrightarrow\!\!\!\!\!\!\!\!&\!\!\!\!\!\!\!  \longrightarrow
&&&&&  \longrightarrow\!\!\!\!\!\!&\!\!\! \longrightarrow\!\!\!\!&\!\!\!\!\!  \longrightarrow\!\!\!\!\!\!&\!\!\!\!\!\!\!\!\! \longrightarrow
\end{array}
\]
\vspace{-0.1cm}
\[
\begin{array}{lllllllllllllllll}
\vspace{-0.3cm}
& \circ &&\!\!\!\!\circ&&\!\!\!\!\circ&&\!\!\!\!\!\circ&\qquad & \circ &&\!\!\!\!\circ&&\!\!\!\!\circ&&\!\!\!\!\!\circ&\\
\vspace{-0.2cm}
\!\!\!o_{A_7}\!\cup\tilde d_3:\!\!\!& \ \ \nwarrow&\nearrow&\!\nwarrow&\nearrow&\!\nwarrow&\nearrow  &&\!\!\!\overset{red}{\longrightarrow}\  \tilde{o d_3}:\!\!\!& \ \ \nwarrow&\nearrow&&&\!\nwarrow&\nearrow   \\
&&\!\!\!\!\circ &\!\!\!\!\!\!\!\!\longleftarrow&\!\!\!\!\bullet&\!\!\!\!\!\!\!\!\longrightarrow&\!\!\!\!\circ&& &&\!\!\!\!\circ &\!\!\!\!\!\!\!\!\!\!\longleftarrow&\!\!\!\!\bullet&\!\!\!\!\!\longrightarrow&\!\!\!\!\circ&&\\
\end{array}
\]
Two more rooted trees are obtained from $\tilde {od_1}$ and $\tilde{o d_2}$ by the 
action of the
left-right involutive diagram automorphism of $\Gamma(A_7)$. 
Therefore, the formula (19) of the principal number of type $A_7$ turns out to be 
\[\begin{array}{lll}
\sigma(A_7)=2\sigma(\tilde{od_1})+2\sigma(\tilde{od_2})+\sigma(\tilde{od_3})
=2\frac{7!}{7\cdot5\cdot3}+2\frac{7!}{7\cdot5\cdot3}+\frac{7!}{7\cdot3\cdot3}
=272.
\\
\end{array}
\]


\noindent 
{\bf II.} \ \  
The principal decomposition of $\Gamma(A_7)$ opposite to (13) 
and the opposite principal orientation $-o_{A_7}:=o_{\Pi_2,\Pi_1}$ are given by 
\[
\begin{array}{llllllllll}
\vspace{-0.2cm}
&\Pi_1:&&&\!\!\!\! \circ&&\!\!\!\!\circ&&\!\!\!\!\circ\\
\vspace{-0.2cm}
\!\!\!\!\!\!-o_{A_7}:&&& \nearrow & \nwarrow&\nearrow&\!\nwarrow&\nearrow&\!\nwarrow  \\
&\Pi_2:&&\!\! \circ&&\!\!\!\!\circ &&\!\!\!\!\circ&&\!\!\!\!\!\circ .\\
\end{array}
\]
There are 14 partial orderings $\tilde d\in  \widetilde{Ord}(\Pi_2)$ and,
accordingly,  
14 associated rooted trees $\tilde{od}$. Seven of them are illustrated as follows: 
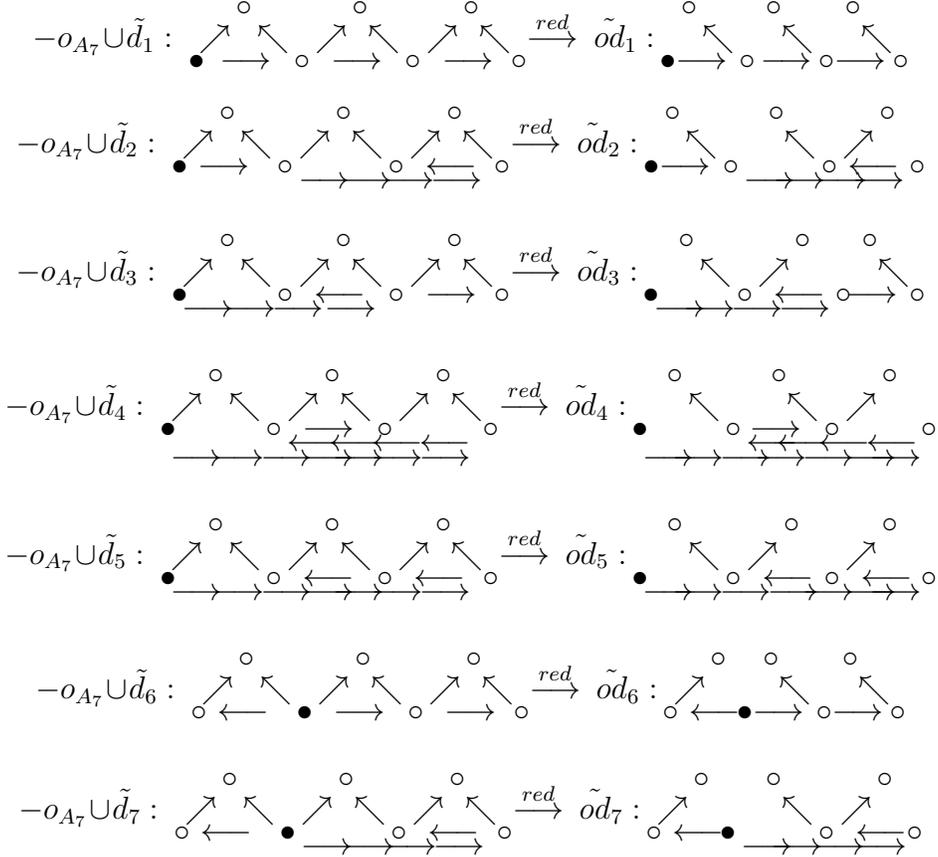
\[
\begin{array}{lllllllllllllllllllll}
\vspace{-0.2cm}
&&\!\!\!\! \circ&&\!\!\!\!\circ&&\!\!\!\!\circ&& 
&&\!\!\!\! \circ&&\!\!\!\!\circ&&\!\!\!\!\circ\\
\vspace{-0.2cm}
\!\!\!\!\!-o_{A_7}\!\cup\! \tilde d_1:& \nearrow & \nwarrow&\nearrow&\!\nwarrow&\nearrow&\!\nwarrow  &\!\!\overset{red}{\longrightarrow}\!\!&\tilde{od}_1:\!\!
&  & \nwarrow&&\!\nwarrow&&\!\nwarrow  \\
&\!\! \bullet&\!\!\!\!\!\!\!\longrightarrow\!\!\!\!&\!\!\!\!\circ &\!\!\!\!\!\!\!\!\longrightarrow\!\!\!\!&\!\!\!\!\circ&\!\!\!\!\!\!\!\!\longrightarrow\!\!\!\!&\!\!\!\!\!\circ &
&\!\! \bullet&\!\!\!\!\!\longrightarrow\!\!\!\!&\!\!\!\!\circ &\!\!\!\!\!\longrightarrow\!\!\!\!&\!\!\!\!\circ&\!\!\!\!\!\!\longrightarrow\!\!\!\!&\!\!\!\!\!\circ \\
\end{array}
\vspace{-0.1cm}
\]
\[
\begin{array}{lllllllllllllllllllll}
\vspace{-0.2cm}
&&\!\!\!\! \circ&&\!\!\!\!\circ&&\!\!\!\!\circ&& 
&&\!\!\!\! \circ&&\!\!\!\!\circ&&\!\!\!\!\circ\\
\vspace{-0.2cm}
\!\!\!\!\!-o_{A_7}\!\cup\! \tilde d_2:& \nearrow & \nwarrow&\nearrow&\!\nwarrow&\nearrow&\!\nwarrow  &\!\!\overset{red}{\longrightarrow}\!\!&\tilde{od}_2:\!\!
& & \nwarrow&&\!\nwarrow&\nearrow&  \\
&\!\! \bullet&\!\!\!\!\!\!\!\!\longrightarrow\!\!\!\!&\!\!\!\!\circ &\!\!\!\!\!\!\!\!\!\!\!\!&\!\!\!\!\circ&\!\!\!\!\!\!\!\!\longleftarrow\!\!\!\!&\!\!\!\!\!\circ &
&\!\! \bullet&\!\!\!\!\!\longrightarrow\!\!\!\!&\!\!\!\!\circ &\!\!\!\!\!\!\!\!\!\!\!\!&\!\!\!\!\circ&\!\!\!\!\!\!\!\!\!\!\longleftarrow\!\!\!\!&\!\!\!\!\circ
\vspace{-0.3cm}
 \\
&&&  \longrightarrow\!\!\!\!&\!\!\! \longrightarrow\!\!\!\!&\!\!\!\!\!  \longrightarrow\!\!\!\!\!\!\!\!&\!\!\!\!\!\!\!  \longrightarrow
&&&&&  \longrightarrow\!\!\!\!\!\!&\!\!\! \longrightarrow\!\!\!\!&\!\!\!\!\!  \longrightarrow\!\!\!\!\!\!&\!\!\!\!\!\!\!\!\! \longrightarrow
\end{array}
\vspace{-0.1cm}
\]

\[
\begin{array}{lllllllllllllllllllll}
\vspace{-0.2cm}
&&\!\!\!\! \circ&&\!\!\!\!\circ&&\!\!\!\!\circ&& 
&&\!\!\!\! \circ&&\!\!\!\!\circ&&\!\!\!\!\circ\\
\vspace{-0.2cm}
\!\!\!\!\!-o_{A_7}\!\cup\! \tilde d_3:&  \nearrow& \nwarrow&\nearrow&\!\nwarrow&\nearrow&\!\nwarrow  &\!\!\overset{red}{\longrightarrow}\!\!&\tilde{od}_3:\!\!
& & \nwarrow&\nearrow&&&\!\nwarrow  \\
&\!\! \bullet&\!\!\!\!\!\!\!\!\!&\!\!\!\!\circ &\!\!\!\!\!\!\!\!\longleftarrow\!\!\!\!&\!\!\!\!\circ&\!\!\!\!\!\!\!\!\longrightarrow\!\!\!\!&\!\!\!\!\!\circ &
&\!\! \bullet&\!\!\!\!\!\!\!\!\!\!&\!\!\!\!\circ &\!\!\!\!\!\!\!\!\longleftarrow\!\!\!\!&\!\!\!\!\circ&\!\!\!\!\!\!\!\longrightarrow\!\!\!\!&\!\!\!\!\!\circ 
\vspace{-0.3cm}
 \\
&  \longrightarrow\!\!\!\!&\!\!\! \longrightarrow\!\!\!\!&\!\!\!\!\!  \longrightarrow\!\!\!\!\!\!\!\!&\!\!\!\!\!\!  \longrightarrow
&&&&&  \longrightarrow\!\!\!\!\!\!&\!\!\! \longrightarrow\!\!\!\!&\!\!\!\!\!  \longrightarrow\!\!\!\!\!\!&\!\!\!\!\!\!\! \longrightarrow
\end{array}
\vspace{-0.1cm}
\]

\[
\begin{array}{lllllllllllllllllllll}
\vspace{-0.2cm}
&&\!\!\!\! \circ&&\!\!\!\!\circ&&\!\!\!\!\circ&& 
&&\!\!\!\! \circ&&\!\!\!\!\circ&&\!\!\!\!\circ\\
\vspace{-0.2cm}
\!\!\!\!\!-o_{A_7}\!\cup\! \tilde d_4:& \nearrow & \nwarrow&\nearrow&\!\nwarrow&\nearrow&\!\nwarrow  &\!\!\overset{red}{\longrightarrow}\!\!&\tilde{od}_4:\!\!
& & \nwarrow&&\!\nwarrow&\nearrow&  \\

&\!\! \bullet&\!\!\!\!\!\!\!\!\!\!&\!\!\!\!\circ &\!\!\!\!\!\!\!\!\longrightarrow\!\!\!\!&\!\!\!\!\circ&\!\!\!\!\!\!\!\!\!\!&\!\!\!\!\!\circ &
&\!\! \bullet&\!\!\!\!\!\!\!\!\!\!&\!\!\!\!\circ &\!\!\!\!\!\!\!\!\longrightarrow\!\!\!\!&\!\!\!\!\circ&\!\!\!\!\!\!\!\!\!\!\!&\!\!\!\!\!\circ 
\vspace{-0.3cm}
 \\
&&&  \longleftarrow\!\!\!\!&\!\!\! \longleftarrow\!\!\!\!&\!\!\!\!\!  \longleftarrow\!\!\!\!\!\!\!\!\!&\!\!\!\!\!\!\!  \longleftarrow
&&&&&  \longleftarrow\!\!\!\!\!\!&\!\!\! \longleftarrow\!\!\!\!&\!\!\!\!\!  \longleftarrow\!\!\!\!\!\!&\!\!\!\!\!\!\! \longleftarrow
\vspace{-0.3cm}
 \\
&  \longrightarrow\!\!\!\!&\!\!\! \longrightarrow\!\!\!\!&\!\!\!\!\!  \longrightarrow\!\!\!\!\!\!\!\!&\!\!\!\!\!\!\!\!  \longrightarrow\!\!\!\! \longrightarrow\!\!\!&\!\!\!\!\!  \longrightarrow\!\!\!\!\!&\!\!\!\!\!\!\!  \longrightarrow
&&&  \longrightarrow\!\!\!\!\!\!&\!\!\! \longrightarrow\!\!\!\!&\!\!\!\!\!  \longrightarrow\!\!\!\!\!\!&\!\!\!\!\!\!\! \longrightarrow&\!\!\!\!\!\!\!\!  \longrightarrow\!\!\!\!\!\!\!\!\!&\!\!\!\!\!\!\!\!\!\!\!  \longrightarrow\!\!\!\!\!  \longrightarrow
\end{array}
\vspace{-0.1cm}
\]

\[
\begin{array}{lllllllllllllllllllll}
\vspace{-0.2cm}
&&\!\!\!\! \circ&&\!\!\!\!\circ&&\!\!\!\!\circ&& 
&&\!\!\!\! \circ&&\!\!\!\!\circ&&\!\!\!\!\circ\\
\vspace{-0.2cm}
\!\!\!\!\!-o_{A_7}\!\cup\! \tilde d_5:& \nearrow & \nwarrow&\nearrow&\!\nwarrow&\nearrow&\!\nwarrow  &\!\!\overset{red}{\longrightarrow}\!\!&\tilde{od}_5:\!\!
& & \nwarrow&\nearrow&&\nearrow&  \\
&\!\! \bullet&\!\!\!\!\!\!\!\!\!\!\!\!&\!\!\!\!\circ &\!\!\!\!\!\!\!\!\longleftarrow\!\!\!\!&\!\!\!\!\circ&\!\!\!\!\!\!\!\!\longleftarrow\!\!\!\!&\!\!\!\!\!\circ &
&\!\! \bullet&\!\!\!\!\!\!\!\!\!\!\!\!&\!\!\!\!\circ &\!\!\!\!\!\!\!\!\longleftarrow\!\!\!\!&\!\!\!\!\circ&\!\!\!\!\!\!\!\!\longleftarrow\!\!\!\!&\!\!\!\!\!\circ 
\vspace{-0.3cm}
 \\
&  \longrightarrow\!\!\!\!&\!\!\! \longrightarrow\!\!\!\!&\!\!\!\!\!  \longrightarrow\!\!\!\!\!\!\!\!&\!\!\!\!\!\!\!\!  \longrightarrow\!\!\!\! \longrightarrow\!\!\!&\!\!\!\!\!  \longrightarrow\!\!\!\!\!&\!\!\!\!\!\!\!  \longrightarrow
&&&  \longrightarrow\!\!\!\!\!\!&\!\!\! \longrightarrow\!\!\!\!&\!\!\!\!\!  \longrightarrow\!\!\!\!\!\!&\!\!\!\!\!\!\! \longrightarrow&\!\!\!\!\!\!\!\!  \longrightarrow\!\!\!\!\!\!\!\!\!&\!\!\!\!\!\!\!\!\!\!\!  \longrightarrow\!\!\!\!\!  \longrightarrow
\end{array}
\vspace{-0.1cm}
\]

\[
\begin{array}{lllllllllllllllllllll}
\vspace{-0.2cm}
&&\!\!\!\! \circ&&\!\!\!\!\circ&&\!\!\!\!\circ&& 
&&\!\!\!\! \circ&&\!\!\!\!\circ&&\!\!\!\!\circ\\
\vspace{-0.2cm}
\!\!\!\!\!-o_{A_7}\!\cup\! \tilde d_6:& \nearrow & \nwarrow&\nearrow&\!\nwarrow&\nearrow&\!\nwarrow  &\!\!\overset{red}{\longrightarrow}\!\!&\tilde{od}_6:\!\!
& \nearrow &&&\!\nwarrow&&\!\nwarrow  \\
&\!\! \circ&\!\!\!\!\!\!\!\!\longleftarrow\!\!\!\!&\!\!\!\!\bullet &\!\!\!\!\!\!\!\!\longrightarrow\!\!\!\!&\!\!\!\!\circ&\!\!\!\!\!\!\!\!\longrightarrow\!\!\!\!&\!\!\!\!\!\circ &
&\!\! \circ&\!\!\!\!\!\!\!\!\longleftarrow\!\!\!\!&\!\!\!\!\bullet &\!\!\!\!\!\!\longrightarrow\!\!\!\!&\!\!\!\!\circ&\!\!\!\!\!\!\longrightarrow\!\!\!\!&\!\!\!\!\!\circ \\
\end{array}
\vspace{-0.1cm}
\]

\[
\begin{array}{lllllllllllllllllllll}
\vspace{-0.2cm}
&&\!\!\!\! \circ&&\!\!\!\!\circ&&\!\!\!\!\circ&& 
&&\!\!\!\! \circ&&\!\!\!\!\circ&&\!\!\!\!\circ\\
\vspace{-0.2cm}
\!\!\!\!\!-o_{A_7}\!\cup\! \tilde d_7:& \nearrow & \nwarrow&\nearrow&\!\nwarrow&\nearrow&\!\nwarrow  &\!\!\overset{red}{\longrightarrow}\!\!&\tilde{od}_7:\!\!
& \nearrow &&&\!\nwarrow&\nearrow&  \\
&\!\! \circ&\!\!\!\!\!\!\!\!\longleftarrow\!\!\!\!&\!\!\!\!\bullet &\!\!\!\!\!\!\!\!\!\!\!\!&\!\!\!\!\circ&\!\!\!\!\!\!\!\!\longleftarrow\!\!\!\!&\!\!\!\!\!\circ &
&\!\! \circ&\!\!\!\!\!\!\!\!\longleftarrow\!\!\!\!&\!\!\!\!\bullet &\!\!\!\!\!\!\!\!\!\!\!&\!\!\!\!\circ&\!\!\!\!\!\!\!\!\longleftarrow\!\!\!\!&\!\!\!\!\!\circ 
\vspace{-0.3cm}
 \\
&&&  \longrightarrow\!\!\!\!&\!\!\! \longrightarrow\!\!\!\!&\!\!\!\!\!  \longrightarrow\!\!\!\!\!\!\!\!&\!\!\!\!\!\!\!  \longrightarrow
&&&&&  \longrightarrow\!\!\!\!\!\!&\!\!\! \longrightarrow\!\!\!\!&\!\!\!\!\!  \longrightarrow\!\!\!\!\!\!&\!\!\!\!\!\!\!\! \longrightarrow
\end{array}
\]

The remaining seven rooted   trees are obtained from 
the above seven trees by the action of the
left-right involutive diagram automorphism of $\Gamma(A_7)$. 
Therefore, the formula (19) of the principal number of type $A_7$ 
turns out to be 
\[
\begin{array}{llll}
\sigma(A_7)\!\!\!&=&\!\!2\big(\sigma(\tilde{od_1})\!+\!\sigma(\tilde{od_2})\!+\!\sigma(\tilde{od_3})\!+\!
\sigma(\tilde{od_4})\!+\!\sigma(\tilde{od_5})\!+\!\sigma(\tilde{od_6})\!+\!
\sigma(\tilde{od_7})\big)
\vspace{0.2cm}
\\
&=&\!\!2\big(\frac{7!}{7\cdot6\cdot4\cdot2}+
\frac{7!}{7\cdot6\cdot4\cdot3}+
\frac{7!}{7\cdot6\cdot3\cdot2}+
\frac{7!}{7\cdot6\cdot5\cdot3}+
\frac{7!}{7\cdot6\cdot5\cdot3}+
\frac{7!}{7\cdot4\cdot2\cdot2}+
\frac{7!}{7\cdot4\cdot3\cdot2}\big)
\vspace{0.2cm}
\\
&=&272.
\end{array}
\]

\noindent
{\it Note.} Even though the starting diagram $\Gamma(A_7)$ is  linear 
the blocks correspond to non-linear diagrams.

\section{Geometric background}

We recall briefly a theorem [S1\S3], 
which combines the 
principal $\Gamma$-cones with some geometry of real bifurcation set 
in case of $\Gamma$ being a Coxeter graph of finite type. 
For details, one is refereed to [S1].
 
\medskip
Let $W$ be a finite reflection group acting irreducibly on an $\R$-vector space 
$V$ of rank $l$.\ Due to Chevalley Theorem [B,ch.v,5.3], 
the invariants  $S(V^*\!)\!^W$\! is freely generated by some 
homogeneous elements $P_1,\!\cdot\!\cdot\!\cdot,\!P_l$.
Thus, the quotient 
variety $S_W\!\!:=\!\!V\!\!/\!\!/W\!\!=\!\!Spec(S(V^*)\!^W)$  
is a smooth affine scheme over $\R$ of coordinates 
$P_1,\!\cdot\!\cdot\!\cdot,P_l$.
It contains the discriminant divisor $D_W$ defined by a polynomial 
$\Delta_W\!\in\! S(V^*)^W$, which is the square of a basic anti-invariant [B,ch.v, 5.4].
{\it The $\Delta_W$ is a monic polynomial of degree $l$ with respect to  
the coordinate $P_l$ of $S_W$ of the largest degree.}
The integration $\exp(D)$\! of the lowest degree vector field 
$D\!:=\!\frac{\partial}{\partial P_l}$ on 
$S_W$ 

\noindent
(unique up to a constant factor and is called 
the {\it primitive vector field}) defines a unique (up to a scaling factor)
additive group $\G_a$-action 
on $S_W$. 
The quotient $T_W\!:=\!S_W/\!\!/\G_a$ is an $(l-1)$-dimensional 
affine variety (forgetting the coordinate $P_l$). 
The restriction to $D_W$ 
of the projection map $S_W\!\to\! T_W$ is an $l$-fold flat covering, whose ramification 
divisor in $T_W$ is denoted by $B_W$ and called\! the {\it bifurcation divisor}.\ 
The\! $B_W$\! decomposes into the ordinary 
part $B_{W,2}$\! and\! the\! higher\! part\! $B_{W,\ge3}$\! according as the 
ramification index of the covering is equal or larger than 2. 

Depending on $\varepsilon\!\in\! \{\pm\!1\}$, there are real forms 
$T_{W,\R}^\varepsilon$,\! $B_{W,2,\R}^\varepsilon$ and 
$B_{W,\ge3,\R}^\varepsilon$ 
of these schemes. 
On the other hand,
arising from a study of eigenspaces of Coxeter elements, 
there is a distinguished real half axis $AO^\varepsilon\!\simeq\!\R_{>0}$ 
embedded in $T_{W,\R}^\varepsilon\!\setminus\! B_{W,\ge3,\R}^\varepsilon$ 
(see [S1] for details). 
The connected component of 
$T_{W,\R}^\varepsilon\!\setminus\! B_{W,\ge3,\R}^\varepsilon$ containing 
 $AO^\varepsilon$, denoted by $E_W^\varepsilon$, is called the 
{\it central region}, 
which we determine in the following theorem.

Using $P_l$, consider the $l$-valued algebroid function  
$T_W\!\!\leftarrow\!\! D_W\!\!\overset{P_l|_{D_W}}{\to}\!\! \A$ 
(=affine line),
defined by the polynomial equation $\Delta_W\!=\!0$ in $P_l$.
It is ramifying
along $B_{W,\ge3}$, but regular along $B_{W,2}$. 
We can indexify  $l$ branches  
of the function on $T^{\varepsilon}_{W,\R}$  
at the base point $AO^\varepsilon$ by the set $\Pi$ of 
a simple generator system of $W$, and  
denote them by $\{\varphi_{\al,\varepsilon}\}_{\al\in\Pi}$ [S1].  
\vspace{-0.1cm}

\begin{theorem} Let $\Gamma(W)$ be the Coxeter graph on $\Pi$, 
and $E_{\Gamma(W)}$ be its principal cone.\ 
Then the algebraic correspondence 
$b_{W,\varepsilon}\!\!:=\!\!\sum_{\al\in\Pi}\varphi_{\al,\varepsilon}\! \cdot\! v_\al$ 
from $T^{\varepsilon}_{W,\R}$ to $V_\Pi$ induces a semi-algebraic homeomorphism: 
\begin{equation}
\!\!\!\!\!\!\! b_{W,\varepsilon}\ :\ \ \overline{E}_{W}^\varepsilon\ \simeq\ \overline{E}_{\Gamma(W)}
\end{equation}
from the closure of the central region of $W$
to the closure of the principal cone for the Coxeter graph $\Gamma(W)$ 
of $W$on $\Pi$, and a homeomorphism: 
\begin{equation}
\ \ \ \ \ b_{W,\varepsilon}\ : \ \ {E}_{W}^\varepsilon\cap B_{W,2,\R}^\varepsilon\  \simeq\  {E}_{\Gamma(W)}\cap
\big( \cup_{\al\beta\in \Pi} H_{\al\beta}\big).\!
\end{equation}

That\! is:\! the central region $E_{W}^\varepsilon$ 
is a simplicial cone homeomorphic to the principal $\Gamma(W)$-cone $E_{\Gamma(W)}$.
Connected components of $E_{W}^\varepsilon\!\setminus\!  B_{W,2,\R}^\varepsilon$
are in one to one correspondence with 
the set 
of chambers 
contained in the principal $\Gamma(W)$-cone $E_{\Gamma(W)}$.
In particular, the number of connected components of 
$E_{W}^\varepsilon\!\setminus\!  B_{W,2,\R}^\varepsilon$ 
is given by 
the principal number $\sigma(\Gamma(W))$.
\end{theorem}
The theorem (in a detailed form) has several important implications 
in the study of the topology of the configuration space $S_W$ ([S1,3]).

\begin{note}
1. {\it 
The correspondence $b_W$ is independent of the choice of $P_l$ 
(up to the scaling constant on the primitive vector field $D$).}

\medskip
\noindent
{\it Proof.} Since the largest exponent of $W$ is unique, any other 
largest degree coordinate $\tilde P_l$ of $S_W$ with $D\tilde P_l=1$ 
is of the form 
$P_l+Q$ for a polynomial $Q$ of lower degree coordinates. Then, 
$\tilde\varphi_\al\!=\!\varphi_\al\! +\!Q$ 
($\al\!\in\!\!\Pi$), whose second term is independent of 
$\al$. So, the ambiguity of translation by $Q$ is absorbed in the 
equivalence in the definition (1) 
of labeled configuration space $V_\Pi$. 
That is: one has $\tilde b_W\!=\!b_W\!+\!Qv_\Pi\!\equiv\! b_W$ in $V_\Pi$. $\Box$

\medskip
2. The principal $\Gamma(W)$-cone in RHS of (31) depends only on 
the underlying graph structure of the diagram $\Gamma(W)$ and not 
on the labels on the edges.
 The graphs $\Gamma(W)$ (forgetting the labels) 
of types $A_l$, $B_l$, $C_l$, $F_4$, $G_2$, 
$H_3$, $H_4$ and $I_2(p)$ are 
 linear. Hence, the central regions $E_W$ for them are homeomorphic 
to the principal cones of type $A_l$. 
\end{note}

\noindent
Finally in the present paper, 
we compare the concept of $\Gamma$-cones with a somewhat similar concept, 
the {\it Springer cones}, which we explain below.

\noindent
{\bf Definition} ([Ar1]).\
Let $V_W$ be a real vector space with an 
irreducible action of a finite reflection group $W$. 
 The reflection hyperplanes of $W$ divide $V_W$ into chambers.
Let $\{H_\al\}_{\al\in\Pi}$ be the system of the walls of a chamber.
A connected component of $V_W\!\!\setminus\!\cup_{\al\in\Pi}H_\al$ 
is called a
{\it Springer cone}. A Springer cone containing the maximal number 
of chambers (unique up to sign [Sp1]) is called a {\it principal Springer 
cone}. This maximal number is called the {\it Springer number}.
The Springer number has been calculated by 
the authors Solomon, Springer and Arnold ([So],[Sp1], [Ar1]). 

There are some formal similarities between 
the (principal) Springer cones in $V_W$ and the (principal) 
$\Gamma$-cones in $V_{\Pi}$ (see Table below). 
A result similar to 
Theorem 3.2 is proven for Springer cones [Sp1, Prop.3].

\medskip
\noindent\!\!\!
\begin{tabular}{|c|c|c|}
\hline
& Springer cone   & $\Gamma$-cone   \\
\hline 
\!\!The ambient\!\!\!&  $V_W$ with $W$-chambers  & 
 $V_\Pi$ with $A_{l-1}$-chambers \\
vector space\!  &\!\!\! (depending on the group $W$)\!\! &\! (depending on the set $\Pi$)\!\\ 
\hline
The cutting   &  $\{H_\al\}_{\al\in\Pi}$ (indexed by 
 &\!\!\! $\!\{H_{\al\beta}\}_{\overline{\al\beta}\in Edge(\Gamma)}$\ (indexed by\!\\
hyperplanes  & the vertices of $\Gamma(W)$) & the edges of the tree $\Gamma$ )\\
\hline
\end{tabular}

 \medskip
Roughly speaking, the principal Springer cones deal with the {\it generators}
of the Artin groups, whereas the principal $\Gamma$-cones 
deal with the {\it (non-commutative) braid relations} of the 
Artin groups.

The only cases when a $\Gamma$-cone decomposition is 
simultaneously a Springer cone decomposition are listed by the following assertion.
\vspace{-0.1cm}

\begin{assertion}
 For a forest $\Gamma$, the following i)--iii) 
are equivalent.

i) There exists  a finite Coxeter group $W$ and a linear isomorphism: 
$V_\Pi\simeq V_W$ which maps chambers to chambers 
and the $\Gamma$-cones to the Springer-cones.


ii) The smallest number of chambers contained in a $\Gamma$-cone 
is equal to 1, i.e.\ $\inf\{\sigma(o) \mid o\in Or(\Gamma)\}=1$. 

iii) $\Gamma$ is a linear graph of type $A_l$, and $W\!=\!W(A_{l-1})$ for
$l\!>\!1$.
\vspace{-0.1cm}
\end{assertion}

\noindent
{\it Proof.}  i)  $\Rightarrow$ ii): This follows from the definition of the 
Springer cone. 

ii) $\Rightarrow$ iii): if a $\Gamma$-cone 
consists of a single chamber 
$\overline C\!:=\!\{\lambda_{\al_1}\!\le\!\cdots\!\le\! \lambda_{\al_l}\}$, 
then $\Gamma$ is a  linear graph 
$\al_1$-$\al_2$-$\cdots$-$\al_l$ (of type $A_l$) on $\Pi$.

iii) $\Rightarrow$ i): 
If $\Gamma$ is a linear graph $\al_1$-$\al_2$-$\cdots$-$\al_l$, 
then the orientation 
$\al_1\!<\! \al_2\! <\!\!\cdots\!\!<\!\al_l$
on $\Gamma$ corresponds to the $\Gamma$-cone consisting only of 
a single chamber  
$\overline C\!:=\!\{\lambda_{\al_1}\!\le\!\cdots\!\le\! \lambda_{\al_l}\}$ 
of type $A_{l-1}$ in $V_\Pi=V_{A_{l-1}}$. \ $\Box$
\begin{remark} Assertion 6.2 is not true if $\Gamma$ is not a forest
(see \S3 Example), since the argument ii) $\Rightarrow$ iii) fails in general.
\end{remark}

Since Assertion 6.2, 
$\sigma(A_l)\!:=\!\sigma(\Gamma(A_l))$ 
is equal to the Springer number $a_{l-1}$ of type $A_{l-1}$. 
Since the Springer number $a_n$ of type $A_{n}$ 
is given by the generating function:
$
 1\!+\!\sum_{n=1}^\infty \frac{a_{n}}{n!}x^{n}\!=\!\frac{1}{1-\sin(x)}
$
([Sp1,3.]), one has
\begin{equation}
1+\sum_{n=1}^\infty \frac{\sigma(A_{n})}{n!}x^{n}=1+\int_0^x
\frac{1}{1-\sin(x)}dx=\mathrm{tan}(\frac{x}{2}+\frac{\pi}{4}).
\vspace{-0.2cm}
\end{equation}
This formula was found by several authors repeatedly (e.g.\ [St, Exercise 43(c)]). 
Including this case,  Y. Sano [Sa] gave 
the following explicit formula of the principal numbers 
for the three infinite series of Coxeter graphs of types 
$A_n$, $D_n$ and $E_n$:
\[\begin{array}{lll}
\sum_{n=1}^\infty \frac{\sigma(A_{n})}{n!}x^{n}&=&\mathrm{tan}(\frac{x}{2}+\frac{\pi}{4})-1,\\
\sum_{n=3}^\infty \frac{\sigma(D_{n})}{n!}x^{n}&=&(2x-1)\mathrm{tan}(\frac{x}{2}+\frac{\pi}{4})+2-2x^2,\\
\sum_{n=4}^\infty \frac{\sigma(E_{n})}{n!}x^{n}&=&(\frac{1}{2}x^2-2x+3)\mathrm{tan}(\frac{x}{2}+\frac{\pi}{4})-3x^2-x-3.
\end{array}\]

\noindent
{\it Question.}
 By an analogy to \S3 Theorem, consider 
any system of $l$-reflection hyperplanes in $V_W$ 
forming coordinate hyperplanes and ask a question:  
is there a unique (up to a sign) orthant of $V_W$, 
cut out by the hyperplanes,  
which contains the maximal number of chambers. 
The answer is apparently positive for the type $A_l$ and $I_2(p)$ for 
odd $p\in 2\Z_{>0}$, and negative 
for the types $B_l,\ C_l$ and $I_2(p)$ for even 
$p\in 2\Z_{>0}$. 

\end{document}